\documentclass[11pt]{article}
\usepackage{graphicx}
\usepackage[utf8]{inputenc}
\usepackage{xcolor}
\usepackage{mathrsfs}
\usepackage{amsmath,amsthm,amssymb,amscd}
\usepackage{booktabs} 
\usepackage{float}
\usepackage[hyperfootnotes=false, colorlinks, linkcolor={blue}, citecolor={magenta}, filecolor={blue}, urlcolor={blue}, plainpages=false, pdfpagelabels]{hyperref}
\usepackage[overload]{textcase} 
\usepackage{mathtools}
\usepackage[numbers,sort&compress]{natbib} 
\usepackage{etoolbox} 
\apptocmd{\sloppy}{\hbadness 10000\relax}{}{} 
\usepackage[left=2.5cm, right=2.5cm, top=3cm]{geometry}
\allowdisplaybreaks
\usepackage{enumerate}
\allowdisplaybreaks[1]
\usepackage[shortlabels]{enumitem}
\usepackage{bm}
\usepackage{url}
\usepackage[multiple]{footmisc}

\usepackage{colonequals}
\usepackage{longtable} 
\usepackage[ampersand]{easylist}
\usepackage{ltablex,booktabs}
\usepackage[
         ordering=Kac,
         edge/.style=black,
         indefinite-edge={draw=green,fill=white,densely dashed},
         indefinite-edge-ratio=2,
         mark=o,
         root-radius=0.1cm]
         {dynkin-diagrams}
\usepackage{stmaryrd}
\usepackage{tikz-cd} 
\usepackage{ytableau}
\usepackage{makecell}

\newcommand{\SL}{{\rm SL}}
\newcommand{\SO}{{\rm SO}}
\newcommand{\Sp}{{\rm Sp}}

\newcommand{\ord}{\mathrm{ord}}
\newcommand{\Ind}{{\rm Ind}}

\newcommand{\forget}[1]{}
\def\qdots{\mathinner{\mkern1mu\raise0pt\vbox{\kern7pt\hbox{.}}\mkern2mu
\raise3.4pt\hbox{.}\mkern2mu\raise7pt\hbox{.}\mkern1mu}}

\usepackage{pifont}

\usepackage[makeroom]{cancel}

\newcommand{\of}[1]{\left(#1\right)}

\newtheorem{lemma}{Lemma}[section]
\newtheorem{theorem}[lemma]{Theorem}
\newtheorem{corollary}[lemma]{Corollary}
\newtheorem{proposition}[lemma]{Proposition}

\newtheorem{remark}[lemma]{Remark}

\newcommand\blfootnote[1]{%
	\begingroup
	\renewcommand\thefootnote{}\footnote{#1}%
	\addtocounter{footnote}{-1}%
	\endgroup
}
\usepackage[toc,page, title, titletoc]{appendix}
\makeatother
\makeatletter
\newcommand\appendix@section[1]{%
	\refstepcounter{section}%
	\orig@section*{Appendix \@Alph\c@section: #1}%
	\addcontentsline{toc}{section}{Appendix \@Alph\c@section: #1}%
}
\g@addto@macro\appendix{\let\section\appendix@section}
\let\orig@section\section
\makeatother
\makeatletter
\def\thickhline{%
  \noalign{\ifnum0=`}\fi\hrule \@height \thickarrayrulewidth \futurelet
   \reserved@a\@xthickhline}
\def\@xthickhline{\ifx\reserved@a\thickhline
               \vskip\doublerulesep
               \vskip-\thickarrayrulewidth
             \fi
      \ifnum0=`{\fi}}
\makeatother

\title{A Rankin-Selberg Integral for $SU_{2,2}$}
\author{Ning Ma}
\date{}

\begin{document}
\maketitle

\blfootnote{\hspace*{0.2in } Key words and phrases: Adjoint $L$-function, Whittaker functions, Unitary groups.}

\begin{abstract}
We construct a Rankin–Selberg integral for the non-split quasi-split unitary group $SU_{2,2}$ and prove that it is nontrivial and Eulerian. We then compute the local integral at the inert places where all data are unramified, and show that it represents the local adjoint L-function of $SU_{2,2}$.
\end{abstract}

\section{Introduction}

In this paper, we construct a Rankin–Selberg integral for the non-split quasi-split unitary group $SU_{2,2}$ by embedding it into the split exceptional group $F_4$, via the isomorphism $SU_{2,2}\cong \operatorname{Spin}^*(6)$. Our construction is based on a modification of the approach of Bump–Ginzburg \cite{B-G1998}. More precisely, by restricting a cuspidal automorphic form on $U_{2,2}$ to $SU_{2,2}$ and integrating it against an Eisenstein series on $F_4$ and a theta function associated with a suitable unipotent subgroup, we obtain an Eulerian integral representation, as stated in Theorem~\ref{main theorem-1}. Furthermore, at an unramified local place, we compute the corresponding local integral and show that it represents the local adjoint $L$-function of $SU_{2,2}$. This is the content of Theorem~\ref{mainTheorem-2}.

\section{The Groups}
In this section we will fix some notation and study a certain double coset space, which will be used in unfolding of the global integral. 

Let $K$ be a number field and let $\mathbb{A}$ be its ring of adeles. In this paper, for any linear algebraic group $G$, we set $[G]=G(K)\backslash G(\mathbb{A})$. Let $E/K$ be a quadratic extension, $R$ be any $K$-algebra and let 
$$J_4=\begin{pmatrix}
 & & & 1\\
 & & 1 & \\
 & 1 & & \\
 1 & & &
\end{pmatrix}.$$ 
We define the special unitary group scheme over $K$ by
$$SU_{2,2}(R)=\{g\in \text{Res}_{E/K}SL_4(R)=SL_4(R\otimes_KE)| gJ_4\overline{g}^t=J_4 \},$$
where $^-$ denotes the conjugation by the nontrivial element of the Galois group $\Gamma_{E/K}$, it is well-known that the $K$-rational points $SU_{2,2}(K)$ is a nonsplit quasisplit unitary group. In the following we realize $E=K(\sqrt{\rho})$ for some $\rho\in K$.

Now, let $F_4$ denote the unique split connected simply connected simple exceptional algebraic group of type $F_4$, we label the positive simple roots of Lie algebra $\frak{f}_4$ as follow:
$$\dynkin[label,label macro/.code={\alpha_{#1}}, edge length=2cm]F4$$
If $\alpha=\pm(n_1\alpha_1+n_2\alpha_2+n_3\alpha_3+n_4\alpha_4)$ with $n_i\in \mathbb{Z}_{\ge 0}$, we abbreviate $\alpha$ by $(n_1n_2n_3n_4)$. We fix a Chevalley basis $\{X_{\alpha}:\alpha\in\Phi\}\cup \{H_{\alpha_i}:i=1,2,3,4\}$ for the Lie algebra $\mathfrak{f}_4$ as follows. First consider a Chevalley basis $\{X'_{\alpha}:\alpha\in\Phi\}\cup \{H'_{\alpha_i}:i=1,2,3,4\}$ with structure constants as in \cite{GZ1988} and then let $H_{\alpha_i}=H'_{\alpha_i}$ for $i=1,2,3,4$ and $X_{\alpha}=X'_{\alpha}$ for $\alpha\neq (1242)$, $X_{1242}=-X'_{1242}$. Finally, we shall denote the structure constants for each pair of roots by $N(\alpha, \beta)$. With the above choice of the above positive roots, we have the standard Borel subgroup $B_{F_4}$ and the maximal split $T_{F_4}$ and the standard parabolic subgroup means it contains $B_{F_4}$. For each root $\alpha$, we denote by $x_{\alpha}: \mathbb{G}_a\to F_4$ the one parameter unipotent subgroup $U_{\alpha}$ or $x_{\alpha}(r)$ corresponding to $\alpha$, and by $h_{\alpha}$ the associated cocharacter, respectively. We shall denote $X_{\pm n_1n_2n_3n_4}$ (resp.$x_{\pm n_1n_2n_3n_4}(r)$) for $X_{\alpha}$ (resp.$x_{\alpha}(r)$). Let $W(F_4)$ be the Weyl group of $F_4$ and $w_i$ be the simple reflection associated to the simple root $\alpha_i$. We shall write $w[i_1i_2...i_r]$ for $w_{i_1}w_{i_2}...w_{i_r}$. In the following we only consider the minimal expression of $w$ for if we realize a Weyl element as a group element, its minimal expression will independent of the choice of the representatives of each simple reflection in $F_4$. Such a representative can be given by for each simple root $\alpha$, let
$\phi_{\alpha}:\SL_2\to F_4$ such that
$$\phi_{\alpha}(\begin{pmatrix}
    1 & 0\\
    t & 1
\end{pmatrix})=x_{-\alpha}(t)\text{, }\phi_{\alpha}(\begin{pmatrix}
    1 & t\\
    0 & 1
\end{pmatrix})=x_{\alpha}(t),$$
then set $\dot{w}_{\alpha}=\phi_{\alpha}(\begin{pmatrix}
    0 & 1\\
    -1 & 0
\end{pmatrix})=x_{\alpha}(1)x_{-\alpha}(-1)x_{\alpha}(1).$

Finally, for the global integral we also need certain subgroups of $F_4$. First, let $P$ be the standard maximal parabolic subgroup associated to the simple roots $\alpha_1, \alpha_3$ and $\alpha_4$, the semisimple part of its Levi subgroup is isomorphic to $SL_2\times SL_3$ and the unipotent radical $R$ is of 20 dimensional, its corresponding roots are those with $n_2\ge 1$. Secondly, let $Q$ be the standard maximal parabolic subgroup associated to the simple roots $\alpha_1, \alpha_2$ and $\alpha_3$, the semisimple part of its Levi subgroup is isomorphic to $Spin_7$ and the unipotent radical $V$ is of 15 dimensional, its corresponding roots are those with $n_4\ge 1$. 
The Levi subgroup of $Q$ which is $GSpin_7$ acts on the Lie algebra $\frak{v}$ of $V$ by adjoint action. On the other hand, as a $GSpin_7$-module, $\frak{v}$ can be decomposed as a direct sum of two subspaces i.e. one is the 7 dimensional subspace  $[\frak{v},\frak{v}]$ on which the semisimple part of $Q$, that is $Spin_7$, acts as standard representation, the other is the 8 dimensional complementary subspace $W$ of $[\frak{v},\frak{v}]$ in $\frak{v}$ which maps isomorphically onto $\frak{v}/[\frak{v},\frak{v}]$ on which $Spin_7$ acts as spin representation. Now, consider the ordered Chevalley basis of $[\frak{v},\frak{v}]$ which is $X_{0122},X_{1122},X_{1222},X_{1232},X_{1242},X_{1342}, X_{2342}$ and let $v_{\rho}=(0,0,1,0,\rho,0,0)$ be a vector with respect to the above Chevalley basis. The isotropy subgroup $G_{\rho}=Stab(v_{\rho})$ in $Spin_7$ is a quasi-split $Spin^*_6$ and thus it's isomorphic as an algebraic group to $SU_{2,2}$ and the simple roots of $G_{\rho}$ are $(1000),(0100),(0120)$, which generates a root system of $A_3$ type. We shall denote the standard Borel subgroup of $SU_{2,2}$ by $B$ and its maximal torus as $T$. With the above embedding, its unipotent radical $N_B$ is given by the following, first, let $M$ be the Levi subgroup of $Q$ and consider its parabolic subgroup which is associated to the simple root $\alpha_3$, then the intersection of $\alpha_3$-unipotent radical with $G_{\rho}$ is our $N_B$. Let $T_d\subset T$ be the maximal split $K$-torus which is contained in a maximal $K$-torus of $SU_{2,2}$, then the relative root system with respect to $T_d$ is of $C_2$-type. Finally, the outer automorphism induced by $^-$ fixes $(1000)$, permutes $(0100)$ and $(0120)$. 

Next, in order to check various assertions, we will consider the matrix representation of general elements of Lie algebras $\mathfrak{so}(7)$ and $\mathfrak{su}(2,2)$ with respect to the above Chevalley basis. Let $[\frak{v},\frak{v}]^{\vee}$ be linear dual space of $[\frak{v},\frak{v}]$ and the dual basis is given by the following linear functional $l_{0122},l_{1122},l_{1222},l_{1232},l_{1242},l_{1342}, l_{2342}$ and therefore in terms of dual space, the Lie algebra $\mathfrak{su}_{2,2}$ is just the stabilizer of the vector $\rho l_{1222}+l_{1242}$ in $\mathfrak{so}_7$. Let $\mathfrak{so}_7$ be its Lie algebra acting on $[\frak{v},\frak{v}]^{\vee}$ by adjoint action, i.e. $ad^{\vee}:\mathfrak{so}_7\to \mathfrak{gl}_7([\frak{v},\frak{v}]^{\vee})$ given by
$$(ad^{\vee}_{X_\alpha}.l_{\beta})(Y) =l_{\beta}(-[X_{\alpha},Y]),$$
for any roots $\alpha$ and $\beta$. Then the matrix representations are the following,

$$\mathfrak{so}(7): \begin{pmatrix}
-T_1 & a & d & u & -y & z & 0\\
a_1 & -T_2+T_1 & -b & -h & x & 0 & -z\\
d_1 & -b_1 & -2T_3+T_2 & -c & 0 & -x & y\\
2u_1 & -2h_1 & -2c_1 & 0 & 2c & 2h & -2u\\
-y_1 & x_1 & 0 & c_1 & 2T_3-T_2 & b & -d\\
z_1 & 0 & -x_1 & h_1 & b_1 & T_2-T_1 & -a\\
0 & -z_1 & y_1 & -u_1 & -d_1 &-a_1 & T_1
\end{pmatrix},$$

$$\mathfrak{spin}^*(6)\cong \mathfrak{su}(2,2): \begin{pmatrix}
-T_1 & a & d & u & -\rho d & z & 0\\
a_1 & -2T_3+T_1 & -b & -h & \rho b & 0 & -z\\
\rho y_1 & -\rho x_1 & 0 & -\rho c_1 & 0 & -\rho b & \rho d\\
2u_1 & -2h_1 & -2c_1 & 0 & 2\rho c_1 & 2h & -2u\\
-y_1 & x_1 & 0 & c_1 & 0 & b & -d\\
z_1 & 0 & -x_1 & h_1 & \rho x_1 & 2T_3-T_1 & -a\\
0 & -z_1 & y_1 & -u_1 & -\rho y_1 &-a_1 & T_1
\end{pmatrix}.$$
Therefore from the above matrix representation, the root subgroups in the standard maximal unipotent subgroup of $SU_{2,2}$ are $x_{1000}(t), x_{0110}(t), x_{1110}(t), x_{1220}(t), x_{1100}(t)x_{1120}(\rho t)$ and $x_{0100}(t)x_{0120}(\rho t)$. Also, by the above matrix representation of $\mathfrak{so}(7)$, the associated symmetric bilinear form is given by
$$J_7=\begin{pmatrix}
 & & & & & & 1\\
 & & & & & 1 &\\
 & & & & 1 & &\\
 & & &\dfrac{1}{2} & & &\\
 & & 1 & & & &\\
 & 1 & & & & &\\
 1 & & & & & &\\
\end{pmatrix}.$$

Finally, for the unfolding of global integral, we also need the following double coset decomposition.

\begin{lemma}\label{lemma2.1}
If $\rho$ is not a square. Then the space $P(K)\backslash F_4(K)/SU_{2,2}(K)\cdot V(K)$ is finite. A set of representatives is given by the following: 
\begin{enumerate}
 \item[(1)] $e, w[2], w[234], w[2342], w[23412], w[231234], w[23123412], w[2342312342], w[23423123412],$
 \item[(2)] $w[234231234], w_0=w[2312343231234], w_0w[2].$

\end{enumerate} 
 \end{lemma}

\begin{proof} 
In the proof we only consider $K$-rational points, so we suppess the "$(K)$". First, we know that a set of representatives of the space $P\backslash F_4/Q$ is given by
$$[P\backslash F_4/Q]=\{e, w[234], w[231234], w[234231234], w_0\},$$ 
see \cite{B-G1998}, p.122. Noticed that each of them is the unique minimal expression in the corresponding double coset. For each $w\in [P\backslash F_4/Q]$, we denote $S_w$ the set of representatives of the double cosets of $w^{-1}Pw\cap GSpin_7\backslash GSpin_7/SU_{2,2}$. Then the set $\{w\gamma:w\in [P\backslash F_4/Q] \text{ and }\gamma\in S_w\}$ is a set of representatives of $P\backslash F_4/SU_{2,2}\cdot V$.

If we consider the standard representation of $GSpin_7$ as left multiplication on column vector space, Let $T$ be the maximal split torus of $GSpin_7$ and the quotient space $Spin_7/
SU_{2,2}$ can be viewed as the orbit space of $v_{\rho}^t$, which consists of all vectors of norm $2\rho$, where the norm is given by $J_7$, i.e. $v_{\rho}J_7v_{\rho}^t$. Now, for each $w$ in $[P\backslash F_4/Q]$, we consider the group $w^{-1}Pw\cap GSpin_7$ acts on $Spin_7/SU_{2,2}$. It is not hard to show for each $w\in [P\backslash F_4/Q]$ we have $wx_{\alpha}w^{-1}\in P$ for $\alpha$ being one of $(1000)$, $(0100)$ and $(0010)$, therefore these groups are all standard parabolic subgroups of $GSpin_7$ and by the same way as in \cite{B-G1998}, p.124, we have the following: 
\begin{enumerate}
  \item For $e$ and $w_0$, it is $(GL_2\cdot GSpin_3)U_1,$
  \item For $w[234]$ and $w[234231234]$, it is $(T\cdot GSpin_3)U_2,$ 
  \item For $w[231234]$, it is $(T\cdot GL_2)U_3.$ Here $GL_2$ is the one generated by $x_{0100}$,
\end{enumerate}
where $U_1$, $U_2$ and $U_3$ are the unipotent subgroups for the corresponding parabolic subgroups. 

Now, let $v$ be a vector in the $Spin_7$-orbit of $v_{\rho}$. For i), we can write $v=(v_1, v_2, v_3)$ with $v_2$ is of 3 dimensional. Either $v_3=0$ or not gives two orbits. This is because if $v_3=0$, $v_2$ cannot be 0 by the length of $v_{\rho}$. So we can choose an element in $U_1$ such that $v_1$ becomes 0, $GSpin_3$ acts on $v_2$ gives the orbit that contains $v_{\rho}$. If $v_3\neq 0$, we can choose element in $U_1$ such that $v_2$ becomes 0, in this case, it contains an vector $(0, 1, 0, 0, 0, \rho, 0)$. On the other hand, $w[2](\alpha_{1222})=\alpha_{1122}$ and $w[2](\alpha_{1242})=\alpha_{1342}$. So in i) we get two double coset representatives for each, $e, w[2], w_0, w_0w[2]$.

For ii), we write $v=(x_1, x_2, v_3, y_1, y_2)$, here $v_3$ is a 3 dimensional vector. If $y_1=y_2=0$, the same argument gives an orbit contains $v_{\rho}$. If $y_2=0$ and $y_1\neq 0$, we first choose an element in $U_2$ such that $v_3$ becomes 0, then we see $x_2\neq 0$, now choose an element in $U_2$ again making $x_1=0$ we see that this orbit contains a vector $(0, 1, 0, 0, 0, \rho, 0)$. If $y_2\neq 0$ and $y_1=0$, the same argument gives the orbit contains$(1, 0, 0, 0, 0, 0, \rho)$. Again, consider the Weyl group elements $w[2]$ and $w[12]$ acts on the Chevalley basis of $[\frak{v},\frak{v}]$. For ii) we have 3 orbits for each, they are $w[234], w[234]w[2], w[234]w[12], w[234231234], w[234231234]w[2], w[234231234]w[12]$.

For iii), we write $v=(x_1, v_1, x_2, v_2, x_3)$, here $v_1$ and $v_2$ are 2 dimensional vectors. If $x_3=0$ and if $v_2=0$, we see $\dfrac{x_2^2}{2}=2\rho$ by length of $v_{\rho}$, this is a contradict to $\rho$ is not a square. So $v_2\neq 0$, we can choose an element in $U_3$ such that $x_1$ and $x_2$ both become to 0, thus we get an orbit which contains $v_{\rho}$. If $x_3\neq 0$, then we can choose an element in $U_3$ makes $v_1, x_2$ and $v_2$ are all 0, so we get an orbit contains $(1, 0, 0, 0, 0, 0, \rho)$. Therefore, we have two representatives $w[231234], w[231234]w[12]$ and complete the proof.

\end{proof}
\begin{remark}\label{rmk-1}
    If $\rho=d^2$ is a square, it is easy to see, from the proof of iii), there is one more double coset, which comes from the Weyl element $w[231234]$. In fact, we can choose an element $\gamma_d\in GSpin_3$ such that $\gamma_d(v_{\rho})=(0,0,0,1,0,0,0)$. Thus the representative of this double coset can be chosen as $w[231234]\gamma_d$. 
 \end{remark}

\section{The Global Integral}

In this section we define the global integral and show it is Eulerian. Let $\varphi$ be a cusp form of $U_{2,2}(\mathbb{A})$ whose restriction on the $SU_{2,2}(\mathbb{A})$ has a Whittaker model, i.e. the integral
$$W_{\varphi}(g)=\int_{[N_B]}\varphi(ng)\psi(n)dn$$
is not identically zero for some additive character $\psi$ of $K\backslash\mathbb{A}$, where $N_B$ is the unipotent radical of the Borel subgroup of $SU_{2,2}$. In this paper, we extended the additive character $\psi$ to $N_B$ as
$$\psi(x_{1000}(r_1)x_{0100}(r_2)x_{0120}(\rho r_2)n')=\psi(\rho(r_1+r_2)).$$

Let $P$ be the parabolic subgroup of $F_4$ defined above, $\delta_P$ be its modular character. Given $s\in \mathbb{C}$, let $\Ind_{P(\mathbb{A})}^{F_4(\mathbb{A})}\delta_P^s$ denote the space of all smooth functions $f_s(g)$ on $F_4(\mathbb{A})$ satisfying $f_s(hg)=\delta_P^s(h)f_s(g)$ for all $h\in P(\mathbb{A})$ and $g\in F_4(\mathbb{A})$. The associated Eisenstein series is defined by the following
$$E_P(g,f_s)=\sum_{\gamma\in P(K)\backslash F_4(K)}f_s(\gamma g).$$
This Eisenstein series converges for $Re(s)$ large enough and admits a meromorphic continuation to the complex plane $\mathbb{C}$.

 For the 8 dimensional spin representation on the vector space $W\cong\frak{v}/[\frak{v},\frak{v}]$, let $l_{\rho}\in [\frak{v},\frak{v}]^{\vee} $ (the dual vector space) corresponding to the vector $(0,0,\rho,0,1,0,0)$. Then we get a linear space $W$ with a non-degenerate bilinear skew-symmetric form $\langle X_1, X_2\rangle_{\rho}=l_{\rho}([X_1,X_2])$, where $X_1, X_2$ are vectors in $W$. With this bilinear form we have a symplectic group $Sp_8^{\rho}$, for convenience we write down the matrix of this bilinear form
 $$J_8=\begin{pmatrix}
 & & & & & & & -2\rho\\
 & & & & & & 2 &\\
 & & & & & 2 & &\\
 & & & & 2\rho & & &\\
 & & & -2\rho & & & &\\
 & & -2 & & & & &\\
 & -2 & & & & & &\\
2\rho & & & & & & &\\
\end{pmatrix}$$
with respect to the following basis of $W$,\\

\vspace{1ex}

\noindent
\makebox[0pt][l]{$(\star)$}\makebox[\textwidth][c]{\{$X_{0001}, X_{0011}, X_{0121}, X_{0111}, X_{1111}, X_{1121}, X_{1231}, X_{1221}\},$}\\
\vspace{1ex}

\noindent where the coefficients are based on the structure constant of $F_4$, $N(0001,1221)=-2$, $N(0011,1231)=2$, $N(0121, 1121)=2$ and $N(0111,1111)=2$, see \cite{GZ1988}. We remark that this is not a symplectic basis.

Let $H_9=W\oplus K$ denote the Heisenberg group with nine variables and the multiplication defined by
$$(w_1,t_1)(w_2,t_2)=(w_1+w_2,t_1+t_2+\frac{1}{2}\langle w_1,w_2\rangle_{\rho}).$$
More explicitly, if we write the element in $H_9$ as $(x_1,x_2,x_3,x_4,y_1,y_2,y_3,y_4,z)$ and product given by
$$(x_1,x_2,x_3,x_4,y_1,y_2,y_3,y_4,z)(x'_1,x'_2,x'_3,x'_4,y'_1,y'_2,y'_3,y'_4,z')
$$
$$=(x_1+x'_1,...,x_4+x'_4,y_1+y'_1,...,y_4+y'_4,z+z'-\rho x_1y'_4+ x_2y_3'+ x_3y_2'+\rho x_4y'_1-\rho y_1x'_4-y_2x_3'-y_3x_2'+\rho y_4x'_1).$$
We shall also write $(\underline{x}, \underline{y},z)$ with $\underline{x}=(x_1,x_2,x_3,x_4)$ and $\underline{y}=(y_1,y_2,y_3,y_4)$. Now, for $u\in V$ we have
$$u=x_{0001}(x_1)x_{0011}(x_2)x_{0121}(x_3)x_{0111}(x_4)x_{1111}(y_1)x_{1121}(y_2)x_{1231}(y_3)x_{1221}(y_4)x_{1222}(r_1)x_{1242}(r_2)u',$$
where $u'$ is a product of all other root subgroups in $V$ in any order. Then we define a map $\tau : V\to H_9$, the Heisenberg group, as
$$\tau(u)=(\underline{x}, \underline{0},0)(\underline{0}, \underline{y},0)(\underline{0}, \underline{0},\rho r_1+r_2)$$
and it is easy to see this map is a group homomorphism. We will let $W^+(\mathbb{A})$ be the maximal isotropic subspace of $W$ that generated by the elements of the basis corresponding roots $(0001), (0011), (0121), (0111)$ and let $W^-(\mathbb{A})$ be the maximal isotropic subspace that generated by the elements of the Chevalley basis corresponding roots $(1111), (1121), (1231), (1221)$.

Let $\psi$ be a non-trivial additive character of $K\backslash\mathbb{A}$ and $\mathcal{S}(W^+(\mathbb{A}))$ be the space of Schwartz-Bruhat functions. The Weil representation $\omega_{\psi}$ of $H_9(\mathbb{A}) \rtimes\widetilde{Sp_8}(\mathbb{A})$ can be realized in the Schr\"{o}dinger model acting on the space $\mathcal{S}(W^+(\mathbb{A}))$. Here are the formulas for the action.
\begin{align*}
   &\omega_{\psi}((\underline{x},0,0)).\phi(\xi)=\phi(\xi+\underline{x}),\\
   &\omega_{\psi}((0,\underline{y},0)).\phi(\xi)=\psi(\langle\xi,\underline{y}\rangle)\phi(\xi),\\
   &\omega_{\psi}((0,0,z)).\phi(\xi)=\psi(z)\phi(\xi),\\
   &\omega_{\psi}(\begin{pmatrix}
       g & 0\\
       0 & g^*\end{pmatrix}, \epsilon).\phi(\xi)=\epsilon\gamma_{\psi}(\det(g))|\det(g)|^{\frac{1}{2}}\phi(\xi\cdot g),\\
   &\omega_{\psi}(\begin{pmatrix}
       I & y\\
       0 & I\end{pmatrix}, \epsilon).\phi(\xi)=\epsilon\psi(\frac{1}{2}\langle\xi, \xi\cdot y\rangle)\phi(\xi).
\end{align*}
In the later local computation, we need computing the factor $\epsilon$ explicitly. We recall this action, see\cite{HKW-96}, \S 1. We denote by $m(g)$ a Siegel-Levi element in the 4th formula, the splitting of the metaplectic cover over the unitary group acts
in the Schr\"{o}dinger model by
$$\omega_{\psi}(m(g)).\phi(\xi)=\gamma_{\psi}(N_{E/K}({\det}_E(g)))\lambda({\det}_E(g))|{\det}_{E}(g)|_E^{\frac{1}{2}}\phi(\xi.g)$$
where $\lambda|_{K^{\times}}=\eta$ and $\eta=\eta_{E/K}$ is the
quadratic character of $K^{\times}$ attached to the field extension $E/K$. For the global Weil index $\gamma_{\psi}$ associated to $\psi$, we only need some facts about its local part and we shall recall these now. For non-Archimedean local field $F$, the local Weil index is defined as follow, first if $a\in F$ let
$$\Gamma_{\psi}(a)=\lim_{n\to \infty}\int_{p^{-n}\mathcal{O}_F}\psi(ax^2)dx,$$
then define
$$\gamma_{\psi}(a)=\frac{\Gamma_{\psi}(a)}{|\Gamma_{\psi}(a)|}.$$
If $\psi$ is an unramified character then it is a well known result that $\gamma_{\psi}(1)=\gamma_{\psi}(a^2)=1$, if the characteristic of the residue field of $F$ is odd, see \cite{Chai2007ANO} or \cite{KudlaNotes}, p.17 for example. The associated theta function is
$$\Tilde{\theta}_{\phi}(hg)=\sum_{\xi\in W^+(K)}\omega_{\psi}(hg).\phi(\xi),$$ where $\xi=(\xi_1,\xi_2,\xi_3,\xi_4)\in K^4$.

Our $SU_{2,2}$ is the stabilizer of $l_{\rho}$ and hence also embedded in to $Sp_8^{\rho}$. On the other hand, the Lie algebra $\mathfrak{so}(7)$ acting on the 8 dimensional subspace $W$ preserves a symmetric bilinear form and therefore the preimage of $SU_{2,2}$ in the metaplectic covering $\widetilde{Sp_8^{\rho}}$ is split by \cite{Kazhdan-1977}. We identify $g\in SU_{2,2}(K)$ with its image $(g, \epsilon)$ in $\widetilde{Sp_8^{\rho}}(K)$. Similarly, $SU_{2,2}(\mathbb{A})$ maps into $\widetilde{Sp_8^{\rho}}(\mathbb{A})$. We can realize each element of $\mathfrak{su}_{2,2}$ as a matrix by the adjoint action $ad^{\vee}: \mathfrak{su}_{2,2}\to\mathfrak{gl}_8(W^\vee)$ given by 
$(ad^{\vee}_{X_\alpha}.\epsilon_{\beta})(Y) =\epsilon_{\beta}(-[X_{\alpha},Y])$ for $Y\in W$, where $W^\vee$ is given a dual basis of $(\star)$. We write this basis as 
$$\{\epsilon_{0001}, \epsilon_{0011}, \epsilon_{0121}, \epsilon_{0111},\epsilon_{1111}, \epsilon_{1121}, \epsilon_{1231}, \epsilon_{1221}\}$$
and it is again not a symplectic basis.

Then we get the matrix representation of the image of a general element of the Lie algebra of the Borel of $\mathfrak{su}_{2,2}$ is 
$$\begin{pmatrix}
 -T_3 & -\rho c_1 & -\rho b   & h       & u & -\rho d & 0 & z\\
 c_1  & -T_3     &   h      & b       & d & u & -z & 0\\
      &          & T_3-T_1 & -c_1    & 0 & a & u & d\\
      &          & \rho c_1& T_3-T_1& a & 0 &  \rho d & -u\\
      &          &          &         & T_1-T_3 & \rho c_1 & -\rho b & h\\
      &          &          &         &  -c_1   & T_1-T_3   & -h & -b \\
      &          &          &         &         &           & T_3 & -c_1\\
      &          &          &         &         &           & \rho c_1 & T_3\\
\end{pmatrix}.$$
Finally, in order to define and compute the global, We need to define Haar measures on the various quotient spaces. We first recall the Iwasawa decomposition, let $G(\mathbb{A})=A(\mathbb{A})\cdot N(\mathbb{A})\cdot K$, where $A$ is the maximal torus of $G$, $N$ is the unipotent radical of the minimal parabolic subgroup $B$ of $G$ and $K$ is the maximal compact subgroup of $G$, then we have the corresponding decomposition on the Haar measure $dg=\delta_B(t)^{-1}dtdudk$. Secondly, for the integral over $U(K)\backslash U(\mathbb{A})$, where $U$ is a unipotent subgroup, we mean the integral over any measurable fundamental domain with respect to the Haar measure on $V(\mathbb{A})$ and the Haar measure on it can be expressed as the product of additive Haar measure on a suitable number of copies of $\mathbb{A}$. Since the integral is independent of the choice of fundamental domain the integral is well-defined. With the given data above, we define the global integral
\begin{align}
   I(\varphi, \phi, f_s)=\int_{[SU_{2,2}]}\int_{[V]}\varphi(g)E_P(ug,f_s)\Tilde{\theta}_{\phi}(\tau(u)g)dudg.\label{global int} 
\end{align}

\subsection{The Relation With Bump-Ginzburg's Global Integral}\label{section3.1}
Before unfolding the integral, we would like to work out the relation between our integral (\ref{global int}) and the integral in \cite{B-G1998}. First notice that $v_{\rho}$ and $G_{\rho}$ are still defined when $\rho$ is a nonzero
square. In this case, $G_{\rho}$ is isomorphic to $SL_4$. 

In fact, let $\rho=d^2$ for some $d\in K$, We show that $SU_{2,2}(K)$ is the same as $SL_4(K)$. This is because if $g\in SL_4(K)$, write $g=A+dB$ where $$A=\frac{1}{2}(J_4(g^t)^{-1}J_4+g), B=\frac{1}{2d}(g-J_4(g^t)^{-1}J_4)\in M_4(K),$$ 
and if we let Galois conjugation $\overline{g}=A-dB$, we have $$\overline{g}^tJ_4g=(A^t-dB^t)J_4(A+dB)=J_4.$$
On the other hand, we know that $GSpin_7(K)$ acts on $[\frak{v},\frak{v}]$ as the standard representation, thus there is a $GSpin_7(K)-$element $A_d$ such that $v_{\rho}A_d=(0,0,0,1,0,0,0)=l_{BG}$. In fact, by the Remark~\ref{rmk-1}, we can choose $A_d=h_{\alpha_4}(t_0)x_0$, where $x_0\in Spin_3(K)$ generated by the root subgroups attached to roots $\pm\alpha_3$, $t_0\in K^{\times}$ and $x_0$ are only depends on $d$. Since isotropy subgroup of the vector $l_{BG}$ in $Spin_7$ is just $SL_4(K)$, so we have $A_d^{-1}SU_{2,2}A_d=SL_4(K)$, where the right hand side is Bump-Ginzburg's $SL_4(K)$  in \cite{B-G1998}.

Secondly, as in \cite{B-G1998}, we can choose another polarization of the symplectic space $W$ and we shall denote this by $W_1=W_1^+\oplus W_1^-$associated with the bilinear form $\langle X_1, X_2\rangle_{BG}=l_{BG}([X_1,X_2]) $ as follow. First given an ordered basis (non-symplectic) 
$$X_{0011},X_{0111},X_{1111},X_{1231},X_{0001},X_{0121},X_{1121},X_{1221},$$
we denote the resulting symplectic space as $W_1$, then let $W_1^+$ and $W_1^-$ be the maximal isotropic subspaces generated by $X_{0011},X_{0111},X_{1111},X_{1231}$ and $X_{0001},X_{0121},X_{1121},X_{1221}$, respectively. Once again, with this bilinear form $\langle \_,\_\rangle_{BG}$, we can associate a Heisenberg group $H_9'$.

Now, by using this polarization, we can map the unipotent radical $V$ to the Heisenberg group $H_9'$ as follow, let $u\in V$, we write
$$u=x_{0011}(x_1)x_{0111}(x_2)x_{1111}(x_3)x_{1231}(x_4)x_{0001}(y_1)x_{0121}(y_2)x_{1121}(y_3)x_{1221}(y_4)x_{1232}(r)u',$$
where $u'$ is a product of all other root subgroups associated roots $(0122),(1122),(1222)$, $(1242)$, $(1342)$, $(2342)$. Then we define a map $\tau_1 : V\to H_9'$ as
$$\tau_1(u)=(\underline{x}, \underline{0},0)(\underline{0}, \underline{y},0)(\underline{0}, \underline{0},r).$$

Let $Sp^{BG}_8(K)$ be the symplectic group with respect to the bilinear form $\langle \_,\_\rangle_{BG}$ and therefore $SL_4(K)$ is embedded in to $Sp^{BG}_8(K)$ so that it is the Levi part of the standard Siegel parabolic subgroup. Hence we can consider the associated Weil representation $\omega^1_{\psi}$ of $H_9'(\mathbb{A})\rtimes \widetilde{Sp_8^{BG}}(\mathbb{A})$ acts on $\mathcal{S}(W_1^+(\mathbb{A}))$. Now, we can construct an linear isomorphism $\sigma^+$ of $W^+$ to $W_1^+$, then by a well-known theorem  of symplectic space, see \cite{Rao1993}, Lemma 2.1, $\sigma^+$ can be extended as a symplectomorphism $\sigma$ between symplectic spaces $W$ and $W_1$ so that $\sigma|_{W^+}: W^+\to W_1^+$ and $\sigma|_{W^-}: W^-\to W_1^-$. In other words, $\sigma$ is changing the polarization from $W^+\oplus W^-$ to $W_1^+\oplus W_1^-$ and symplectic form $\langle \_,\_ \rangle_{\rho}$ to $\langle \_,\_ \rangle_{BG}$. Furthermore, it induces an isomorphism between the Schwartz spaces $\sigma^*: \mathcal{S}(W_1^+(\mathbb{A}))\to \mathcal{S}(W^+(\mathbb{A}))$ by $f\mapsto f\circ \sigma^+$ and also an automorphism on the Heisenberg groups $\sigma_1: H_9\to H_9^1$ that is given by $(x,t)\to (\sigma(x),t)$. On the other hand, the linear map $\sigma$ induces an isomorphism between the isometry groups $Sp_8^{\rho}$ and $Sp_8^{BG}$ we shall denote this isomorphism as $\alpha$. Now we have the following commutative diagram
\[\begin{tikzcd}
W \arrow[r, "\sigma"] \arrow[d, "g"] & W_1 \arrow[d, "\alpha(g)"] \\
W \arrow[r, "\sigma"] & W_1
\end{tikzcd}\]
where the vertical arrow stands for the action of $Sp_8^{\rho}$ and $Sp_8^{BG}$. Then this will produce a well-defined group isomorphism $\eta: H_9\rtimes Sp_8^{\rho}\cong H_9^1\rtimes Sp_8^{BG}$ by $\eta(h,g)=(\sigma_1(h), \alpha(g))$, here as usual the conjugation of the semi-direct product $g^{-1}(x,t)g$ is just $(xg, t)$ for any $h=(x,t)\in H_9$.

Now, for each place $v$ we have local field $K_v$, if we realize the nontrivial 2-fold covering group as $Sp_8^{\rho}(K_v)\times \{\pm 1\}$  with the multiplication law as
$$(g_1, \epsilon_1)(g_2, \epsilon_2)=(g_1g_2, c_R(g_1, g_2)\epsilon_1\epsilon_2),$$
where $c_R$ denotes the Ranga Rao's cocycle, see \cite{Rao1993}. Then we can extend the isomorphism $\alpha$ as an isomorphism between covering groups, we have the commutative diagram
\begin{equation}\label{**}
\begin{tikzcd}
\widetilde{Sp_8^{\rho}(K_v)} \arrow[r, "\widetilde{\alpha}"] \arrow[d, "p"] & \widetilde{Sp_8^{BG}(K_v)} \arrow[d, "p_1"] \\
Sp_8^{\rho}(K_v) \arrow[r, "\alpha"] & Sp_8^{BG}(K_v)
\end{tikzcd}
\end{equation}
where $p$ and $p_1$ are covering maps. The reason for this is the following. We know that the uniqueness of the nontrivial 2-fold covering group of symplectic groups see \cite{KudlaNotes}, p.5. Now consider group $\widetilde{Sp_8^{BG}}(K_v)$, we have two covering maps $\alpha\circ p$ and $p_1$ of $Sp_8^{BG}(K_v)$, so by the uniqueness we have the isomorphism $\widetilde{\alpha}$ and therefore the commutative diagram. By using the restricted tensor product, see \cite{Howe1979}, p.276, one can get the adelic version of the above commutative diagram by replacing $K_v$ with $\mathbb{A}$. 

Furthermore, this also gives rise to an isomorphism $\widetilde{\eta}: H_9(\mathbb{A})\rtimes \widetilde{Sp_8^{\rho}(\mathbb{A})}\to H_9^1(\mathbb{A})\rtimes \widetilde{Sp_8^{BG}(\mathbb{A})}$. Moreover, according to the formulas of Weil's representation above, this isomorphism induces the following commutative diagram at the level of Weil's representations
\[\begin{tikzcd}
\mathcal{S}(W_1^+(\mathbb{A})) \arrow[r, "\sigma^*"] \arrow[d, "\omega_{\psi}^1(\widetilde{\eta}(hg))"] & \mathcal{S}(W^+(\mathbb{A})) \arrow[d, "\omega_{\psi}(hg)"] \\
\mathcal{S}(W_1^+(\mathbb{A})) \arrow[r, "\sigma^*"] & \mathcal{S}(W^+(\mathbb{A}))
\end{tikzcd}\]
that is $$\omega_{\psi}^1(\widetilde{\eta}(hg)).\phi_1(\sigma^{-1}(\xi))=\sigma^*(\omega_{\psi}^1(\widetilde{\eta}(hg)).\phi_1(\xi_1))=\omega_{\psi}(hg).\sigma^*(\phi)(\xi).$$
Here we write $\xi=\sigma^+(\xi_1)$ as defined above. Therefore the relation between the theta-function can be derived as follow
\begin{align}
    \Tilde{\theta}_{\phi}(hg)&=\sum_{\xi\in W^+(K)}\omega_{\psi}(hg).\phi(\xi)\nonumber\\
    &=\sum_{\xi_1\in W_1^+(K)}\omega_{\psi}(\Tilde{\eta}^{-1}(h_1g_1)).\phi(\sigma^{-1}(\xi_1))\nonumber\\
    &=\sum_{\xi_1\in W_1^+(K)}\omega^1_{\psi}(h_1g_1).\phi_1(\xi_1)=\Tilde{\theta}_{\phi_1}(h_1g_1)=\Tilde{\theta}^1_{\sigma^{*-1}(\phi)}(\widetilde{\eta}(hg)).\label{Theta}
\end{align}

Now, Let $I_{BG}(\varphi,f_s)$ be the global integral constructed in \cite{B-G1998}, then we have
\begin{proposition}
    If $\rho=d^2\in K$. Then 
    $$I(\varphi, \phi, f_s)=c_dI_{BG}(\varphi_1,f_s^1*\phi_1),$$ 
    where $c_d$ is a positive constant, $\varphi_1(g_1)=\varphi(A_dg_1A_d^{-1})$ is a cusp form on $SL_4(\mathbb{A})$, $f_s^1(g)=f_s(gA_d^{-1})$ is a not necessary $K$-finite section in $\Ind_{P(\mathbb{A})}^{F_4(\mathbb{A})}\delta_P^{s}$ and $\phi_1=\sigma_{BG}^{*-1}(\phi)\in \mathcal{S}(W_1^+(\mathbb{A}))$ here $\sigma_{BG}^{*}$ is an intertwining operator between the two Schr\"{o}dinger models.
\end{proposition}
\begin{proof} Let $\sigma_0=Ad(A_d^{-1})|_W$ and it is easy to see $\sigma_0$ is a symplectic isomorphism between $W$ and $W_1$, thus it give rise to a homomorphism $\sigma_{BG}$ between Heisenberg groups, and also a homomorphism $\alpha$ of symplectic groups as explained on the top of page 9. In order to deal with the changing of the variables, we have $\sigma_{BG}(\tau(u))=\tau_1(A_d^{-1}uA_d)$, i.e. the following commutative diagram of homomorphism of groups
\begin{equation*}
\begin{tikzcd}
V \arrow[r, "\tau"] \arrow[d, "A_d"] & H_9 \arrow[d, "\sigma_{BG}"] \\
V \arrow[r, "\tau_1"] & H'_9
\end{tikzcd}\label{**}
\end{equation*}
On the other hand, the map $\sigma_0$ also induces a map $\sigma_{BG}^*$ between the corresponding Schwartz spaces,  Therefore we have a well-defined isomorphism $\widetilde{\eta_{BG}}:H_9\rtimes \widetilde{Sp_8^{\rho}}\to H_9^1\rtimes \widetilde{Sp_8^{BG}}$ given by $\widetilde{\eta_{BG}}(\tau(u)g)=(\sigma_{BG}(\tau(u)),\widetilde{\alpha}(g))$. Since $\sigma_0$ and $\sigma$ differ by an element
of $Sp^1_8(K)$, we have $\tilde\eta=\tilde\gamma\,\widetilde{\eta_{BG}}(\cdot)\,\tilde\gamma^{-1}$ for some $\tilde\gamma\in\widetilde{Sp}{}^{BG}_8(K)$, so
$\sigma^{*-1}_{BG}:=\omega^1_\psi(\tilde\gamma^{-1})\circ\sigma^{*-1}$ intertwines the two Schr\"odinger models. Since $\tilde\theta^1$ is left invariant under $\widetilde{Sp}{}^{BG}_8(K)$, the identity of theta-functions (\ref{Theta}) gives

\begin{align}
    \Tilde{\theta}_{\phi}(\tau(u)g)=\Tilde{\theta}^1_{\sigma_{BG}^{*-1}(\phi)}(\widetilde{\eta_{BG}}(\tau(u)g)).\label{Theta-1}
\end{align}
By applying the above constructions, we have
\begin{align*}
&\int_{[SU_{2,2}]}\int_{[V]}\varphi(g)E_P(ug,f_s)\Tilde{\theta}_{\phi}(\tau(u)g)dudg\\
=&\int_{[SL_4]}\int_{[V]}\varphi(A_dg_1A_d^{-1})E_P(uA_dg_1A_d^{-1},f_s)\Tilde{\theta}_{\phi}(\tau(u)A_dg_1A_d^{-1})dud(A_dg_1A_d^{-1})\\
=&c_d\int_{[SL_4]}\int_{[V]}\varphi_1(g_1)E_P(uA_dg_1A_d^{-1},f_s)\Tilde{\theta}_{\phi}(\sigma_{BG}^{-1}(\tau_1(A_d^{-1}uA_d))A_dg_1A_d^{-1})dudg_1,\\
\end{align*}
where $\varphi_1(g_1)=\varphi(A_dg_1A_d^{-1})$, the positive constant $c_d$ comes from the changing variable of $g$, $\tau(u)=\sigma_{BG}^{-1}(\tau_1(A_d^{-1}uA_d))$ just follows from the definitions of $\sigma_{BG}$, $\tau$ and $\tau_1$. By the definition of $\alpha$, we have $\alpha(g)=A_d^{-1}gA_d$. Now the last integral is
\begin{align*}
&\int_{[SL_4]}\int_{[V]}\varphi_1(g_1)E_P(uA_dg_1A_d^{-1},f_s)\Tilde{\theta}_{\phi}(\sigma_{BG}^{-1}(\tau_1(A_d^{-1}uA_d))\alpha^{-1}(g_1))dudg_1\\
=&\int_{[SL_4]}\int_{[V]}\varphi_1(g_1)E_P(uA_dg_1A_d^{-1},f_s)\Tilde{\theta}_{\phi}(\widetilde{\eta}_{BG}^{-1}(\tau_1(A_d^{-1}uA_d)g_1))dudg_1\\
=&\int_{[SL_4]}\int_{[V]}\varphi_1(g_1)E_P(uA_dg_1A_d^{-1},f_s)\Tilde{\theta}_{\sigma_{BG}^{*-1}(\phi)}^1(\tau_1(A_d^{-1}uA_d)g_1)dudg_1\\
=&\int_{[SL_4]}\int_{[V]}\varphi_1(g_1)E_P(A_d^{-1}uA_dg_1A_d^{-1},f_s)\Tilde{\theta}_{\sigma_{BG}^{*-1}(\phi)}^1(\tau_1(A_d^{-1}uA_d)g_1)dudg_1.
\end{align*}
Here in the first equality, we apply the definition of $\widetilde{\eta_{BG}}$, in the second equality we use the identity (\ref{Theta-1}) and it makes sense because the image of $SL_4$ embedded in $Sp_8^{BG}$ is semisimple part of the Levi of the Siegel parabolic subgroup and it is well-known that the metaplectic cover splits over $SL_4\cdot N_{Siegel} $ for Ranga Rao's cocycle is trivial, see \cite{Gomez_2021}, p.254. 
The last equality follows from the fact that $A_d^{-1}\in F_4(K)$ and the Eisenstein series is an automorphic form.

Now, let $f_s^1(g)=f_s(A_dgA_d^{-1})=\delta_P^s(A_d)f_s(gA_d^{-1})=\delta_P^s(h_{\alpha_4}(t_0))f_s(gA_d^{-1})=f_s(gA_d^{-1})$, let $u\in V=u^+u'$, where $u^+$ is the product of root subgroups attached to roots $(0011)$, $(0111)$, $(1111)$, $(1231)$ and $u'\in U_{BG}$, which is the unipotent subgroup that is denoted $U$ in the global integral in \cite{B-G1998}. After changing the variable of $u$, then the last integral will be
\begin{align*}
&\int_{[SL_4]}\int_{[V]}\varphi_1(g_1)E_P(ug_1,f_s^1)\Tilde{\theta}^1_{\phi_1}(\tau_1(u)g_1)du_1dg_1\\
=&\int_{[SL_4]}\int_{[U_{BG}]}\int_{(K\backslash \mathbb{A})^4}\varphi_1(g_1)E_P(u^+u'g_1,f_s^1)\Tilde{\theta}^1_{\phi_1}(\tau_1(u^+)\tau_1(u')g_1)du^+du'dg_1.
\end{align*}
Here we let $\phi^1=\sigma_{BG}^{*-1}(\phi)$. Then by the definition of the theta function and $E_P(g,f_s)$ is an automorphic form, we have the last integral is
\begin{align*}
&\int_{[SL_4]}\int_{[U_{BG}]}\varphi_1(g_1)\of{\int_{(K\backslash \mathbb{A})^4}E_P(u^+u'g_1,f_s^1)\omega^1_{\psi}(\tau_1(u^+)\tau_1(u')g_1).\phi_1(0)du^+}du'dg_1\\
=&\int_{[SL_4]}\int_{[U_{BG}]}\varphi_1(g_1)\of{\sum_{\gamma\in P(K)\backslash F_4(K)}\int_{W_1^+(\mathbb{A})}f_s^1(\gamma u'g_1u^+)\phi_1(\tau_1(u^+))du^+}\psi_{U_{BG}}(r_{1232})du'dg_1.
\end{align*}
In the first equality, we apply the property of Weil representation and collapsing the $u^+$ integral with the $\xi$ sum. In the last equality, we first push $u^+$ to the most right and then apply the formulas of Weil representation. Then we assume that we are in the domain of convergence for the Eisenstein series and then plugging in the summation formula of $E_P$.  Finally, let $$(f_s^1*\phi_1)(g)=\int_{ W^+_1(\mathbb{A})}f_s^1(gu^+)\phi_1(\tau_1(u^+))du^+,$$ then we have the last integral is
$$\int_{[SL_4]}\int_{[U_{BG}]}\varphi_1(g_1)E_P(u'g_1, f_s^1*\phi_1)\psi_{U_{BG}}(r_{1232})du'dg_1.$$
This is the global integral in \cite{B-G1998} at page 128. 

\end{proof}

\subsection{The Unfolding}

In this section we will assume $\rho\in K$ is not a square. We will show that the global integral (\ref{global int}) is Eulerian. In fact, we have

\begin{theorem}\label{main theorem-1}
    The integral ~\eqref{global int} converges absolutely for all $s\in \mathbb{C}$ except possibly for those values of $s$ for which $E_P(g,f_s)$ has a pole. For $Re(s)$ sufficiently large, 
    $$I(\varphi, \phi, f_s)=\int_{N_0(\mathbb{A}_K)\backslash SU_{2,2}(\mathbb{A}_K)}\int_{V_0(\mathbb{A}_K)}W_{\varphi}(g)f_s(\gamma_0 u_0g)\omega_{\psi}(\tau(u_0)g).\phi(0,0,0,1)du_0dg,$$
    where $N_0$ is a five dimensional unipotent subgroup of $SU_{2,2}$ generated by $x_{1000}(r)$, $x_{0110}(r)$, $x_{1110}(r)$, $x_{1220}(r)$ and $x_{1100}(r)x_{1120}(\rho r)$. $V_0$ is a 13 dimensional subgroup of $V$ omitting the root subgroups $x_{1121}(r)$ and $x_{1231}(r)$ and $\gamma_0=w_0w[2]$ is a Weyl group element. 
\end{theorem}
\begin{proof}
 Assume $Re(s)\gg 0$, We use the definition of moderate growth, for functions on finite central extension of a reductive group $G(\mathbb{A})$, see \cite{moeglin1995spectral}, p.25. This involves fixing a norm on $G(\mathbb{A})$ as in section I.2.2, which depends on an embedding $i_G$ of $G$ into $SL_n$ for some value of $n$. It follows from properties (v) and (vi) in section I.2.2, that the notion of “moderate growth” is independent of the choice of $i_G$. In particular, we can define two norms on $G_{\rho}$, one obtained by composing $i_{F_4}$ with the embedding of $G_{\rho}$ into $F_4$ and the other obtained by composing $i_{Sp_8^{\rho}}$ with the embedding of $G_{\rho}$ into $Sp_8^{\rho}$. 

The Eisenstein series $E_P(g, f_s)$ as a function of $g$ is a moderate growth function on $F_4(\mathbb{A})$ (with bounds that are uniform as $s$ varies in a compact set). The norm on $F_4$ is continuous, and one of two norms on $G_{\rho}$ is its restriction to $G_{\rho}(\mathbb{A})$. Using these facts, the compactness of $V(K)\backslash V(\mathbb{A})$, and property (ii) from section I.2.2., we see that the function $g\mapsto sup\{|E_P (ug,f_s)|: u\in K\}$ is a moderate-growth function on $G_{\rho}(\mathbb{A})$ for any fixed fundamental domain F for $V(K)/V(\mathbb{A})$ (with bounds that are uniform as s varies in a compact set). Similarly, $\widetilde{\theta}_{\phi}$ is a moderate growth function in $Sp_8^{\rho}$, and therefore  $g\mapsto sup\{|\widetilde{\theta}_{\phi}(ug)|: u\in F\}$ is a moderate growth function in $G_{\rho}(\mathbb{A})$.

As a cusp form $\varphi$ is rapidly decreasing, in the sense of \cite{moeglin1995spectral}, section I.2.12, as follows from Lemma I.2.10. It follows readily from the definition of “rapidly decreasing” together with property (vi) in I.2.2 that the product of a rapidly decreasing function and a moderate growth function is
again rapidly decreasing. Thus $\varphi(g)\cdot sup\{|E_P (ug,f_s)\widetilde{\theta}_{\phi}(ug)|u\in K\}$ is rapidly decreasing and hence so is
$$
   \varphi(g)\int_{[V]}E_P(ug,f_s)\Tilde{\theta}_{\phi}(ug)dudg.
$$
Convergence of the integral follows.
 
Now, we start with unfolding the Eisenstein series. To simplify the notations let $G_{\rho}=SU_{2,2}$ and for each $\gamma$ in the list of Lemma \ref{lemma2.1} we shall write $G_{\rho}^{\gamma}=G_{\rho}\cap\gamma^{-1}P\gamma$ and $V^{\gamma}=V\cap \gamma^{-1}P\gamma$.  Similarly, we define $R_{\rho}=G_{\rho}\cdot V$ and $R_{\rho}^{\gamma}=G_{\rho}\cdot V\cap \gamma^{-1}P\gamma$. With this notations, we have 
 $$I=I(\varphi, \phi, f_s)=\sum_{\gamma\in P\backslash F_4/SU_{2,2}V(K)}\int_{G_{\rho}^\gamma(K)\backslash G_{\rho}(\mathbb{A})}\int_{V^\gamma(K)\backslash V(\mathbb{A})}
 \varphi(g)f_s(\gamma ug)\Tilde{\theta}_{\phi}(\tau(u)g)dudg.$$
 For each representative $\gamma$ we shall compute the integral over the corresponding quotient space. We first show that for all $\gamma$ listed in part (1)
of that Lemma \ref{lemma2.1}, the corresponding integral is 0. In fact, one can check that either the root subgroups $U_{1222}$ or $U_{1242}$ is contained in $R_{\rho}^{\gamma}$ for each $\gamma$ that listed in the part of (1). Then by the property of Weil representation, we have
$$\Tilde{\theta}_{\phi}(\tau(u)g)=\Tilde{\theta}_{\phi}(\tau(x_{1222}(r_1)x_{1242}(r_2))\tau(u_1)g)=\psi(\rho r_1+r_2)\Tilde{\theta}_{\phi}(\tau(u_1)g),$$ 
where we write $u\in V$ as $x_{1222}(r_1)x_{1242}(r_2)u_1$ and $u_1$ is an element of the product of all root subgroups of $V$ other than $U_{1222}$ and $U_{1242}$.
In this case, one will end up with an integral 
 $$\int_{K\backslash\mathbb{A}}\psi(r)dr=0 \text{ or }\int_{K\backslash\mathbb{A}}\psi(\rho r)dr=0.$$
Therefore, the contribution on those double cosets are all 0. In fact, for all coset representatives $\gamma$ listed in (1) in Lemma \ref{lemma2.1}, we have $U_{1242}$ is contained in $V^{\gamma}$.

Next, for $w_1=w[234231234]$ and $w_0$, we will end up with an integral of cusp form $\varphi(g)$ over an unipotent radical of $G_{\rho}$. First, we claim that for those two Weyl elements, the unipotent radical $N_B$ of $G_{\rho}$ is contained in $
G_{\rho}^{w_1}(K)$. Indeed, $N_B$ is generated by the following 1-parameter unipotent subgroups
$$x_{1000}(r), \ x_{0110}(r), \ x_{1110}(r), \ x_{1220}(r), \ x_{0100}(r)x_{0120}(\rho r), \ x_{1100}(r)x_{1120}(\rho r).$$
Then the claim follows from the table~\ref{table 1}. We write $g=n_1g_1\in G_{\rho}$ where $n_1\in N_B$ and we denote the subgroup consists of $g_1$ by $B_{-}$ and we also notice that $N_B$ normalizes $U^-$.

\begin{table}
\centering
\begin{tabular}{|c || c | c | c| c| c| c| c| c|}
 \hline
    \makecell{Positive roots relevant to the standard\\ 
    maximal unipotent of $SU_{2,2}$}
         & 1000 & 0110 & 1110 & 1220 & 0100 & 0120 & 1100 & 1120 \\ [0.2ex]\hline
      $w_1$ & 0122 & 1110 & 1232 & 2342 & 1100 & 1120 & 1222 & 1242 \\ \hline 
      $w_0$ & 1000 & 0121 & 1121 & 1242 & 0120 & 0122 & 1120 & 1122 \\ \hline 
  $w_0w[2]$ & 1120 & 0001 & 1121 & 1122 & -0120 & 0122 & 1000 & 1242 \\ \hline 
\end{tabular}
\caption{$w_1$, $w_0$ and $w_0w[2]$ acts on the roots of unipotent radical of $SU_{2,2}$}
\label{table 1}
\end{table}

Notice that we also have $V^{w_1}(K)$ contains 1-parameter unipotent subgroups $x_{1111}$, $x_{1121}$, $x_{1231}$ and $x_{1221}$, therefore we write $u=u_1u'\in V$ where $u_1$ is the product of root subgroups $U_{1111}$, $U_{1121}$, $U_{1231}$ and $U_{1221}$, we write this unipotent subgroup by $U^-$ and $u'$ is the product of other root subgroups of $V$ and denote $V^-$ therefore
\begin{align*}
&\int\limits_{G_{\rho}^{w_1}(K)\backslash G_{\rho}(\mathbb{A})}\int\limits_{V^{w_1}(K)\backslash V(\mathbb{A})}
 \varphi(g)f_s(w_1 ug)\Tilde{\theta}_{\phi}(\tau(u)g)dudg\\
 =&\int\limits_{B_{-}^{w_1}(K)\backslash B_{-}(\mathbb{A})}\int\limits_{V^{w_1}(K)\backslash V(\mathbb{A})}\int\limits_{[N_B]}
 \varphi(n_1g_1)f_s(w_1 un_1g_1)\Tilde{\theta}_{\phi}(\tau(u)n_1g)dn_1dudg\\
 =&\int\limits_{B_{-}^{w_1}(K)\backslash B_{-}(\mathbb{A})}\int\limits_{V^{w_1}(K)\backslash V(\mathbb{A})}\int\limits_{[N_B]}
 \varphi(n_1g_1)f_s(w_1(n_1^{-1}un_1)g_1)\Tilde{\theta}_{\phi}(n_1\tau(n_1^{-1}un_1)g)dn_1dudg\\
 =&\int\limits_{B_{-}^{w_1}(K)\backslash B_{-}(\mathbb{A})}\int\limits_{V^{w_1}(K)\backslash V(\mathbb{A})}\int\limits_{[N_B]}
 \varphi(n_1g_1)f_s(w_1ug_1)\Tilde{\theta}_{\phi}(n_1\tau(u)g)dn_1dudg
 \end{align*}
 Here we use the fact that for the flat section $f_s$ we have $f_s(w_1un_1g_1)=f_s((w_1n_1w_1^{-1})w_1(n_1^{-1}un_1)g_1)=f_s(w_1(n_1^{-1}un_1)g_1)$ and $N_B$ normalize the unipotent subgroup $V$ so we can change the variable $u$. Next, we consider the integral over $U^-$ we have 
 \begin{align*}
 \int\limits_{[U^-]}f_s(w_1u_1u'g_1)\Tilde{\theta}_{\phi}(n_1\tau(u_1u')g)du_1=f_s(w_1u'g_1)\int\limits_{[U^-]}\Tilde{\theta}_{\phi}(n_1\tau(u_1)\tau(u')g)du_1.
\end{align*}
The last equation is because the invariant property of $f_s$ after pushing $u_1$ across $w_1$ and the change the variable. We carry out the $u_1$ integral on $\widetilde{\theta}_{\phi}$,
\begin{align*}
    \int\limits_{(K\backslash\mathbb{A})^4} \Tilde{\theta}_{\phi}(n_1\tau(u_1)\tau(u')g)du_1=&\int\limits_{(K\backslash\mathbb{A})^4} \Tilde{\theta}_{\phi}(\tau(n_1u_1n_1^{-1})n_1\tau(u')g_1)du_1\\
    =&\int\limits_{(K\backslash\mathbb{A})^4} \sum_{\underline{\xi}\in K^4}\omega_{\psi}(\tau(u_1)n_1\tau(u')g_1).\phi(\underline{\xi})du_1\\
    =&\sum_{\underline{\xi}\in K^4} \omega_{\psi}(n_1\tau(u')g_1).\phi(\underline{\xi}) \int\limits_{(K\backslash\mathbb{A})^4}\psi(-2\rho\xi_4 y_1+2\xi_3 y_2+2\xi_2 y_3+2\rho \xi_1 y_4)d\underline{y}\\
    =&\omega_{\psi}(n_1\tau(u')g_1).\phi(0)=\omega_{\psi}(\tau(u')g_1).\phi(0).
\end{align*}
Here we use the embedding of $G_{\rho}$ into $Sp_8^{\rho}$ above and the last two identities of Weil representation. Finally, put all these together, we end up with the integral of cusp form $\varphi$ over $N_B$ and thus 
\begin{align*}
\int_{[N_B]}\varphi(n_1g_1)dn_1=0
\end{align*}
by the cuspidality of $\varphi$.

For $w_0$, we let $N_P$ be the unipotent radical  of parabolic subgroup of $G_{\rho}$ associated to root $1000$. It is generated by
$$x_{0110}(r), \ x_{1110}(r), \ x_{1220}(r), \ x_{0100}(r)x_{0120}(\rho r), \ x_{1100}(r)x_{1120}(\rho r).$$
This is a subgroup of $N_B$ and therefore is contained in $G_{\rho}^{w_0}$. One can check that $V^{w_0}$ contains only 1-parameter unipotent subgroups $U_{1231}$ and $U_{1221}$ thus, we first carry out the integral over these $U_{1231}$ and $U_{1221}$. Write $u=u_1u'\in V$ where $u_1$ is the product of unipotent subgroups $U_{1231}$ and $U_{1221}$ we denote this 2-dimensional unipotent subgroup by $U_2$ and $u'$ is the product of other 1-parameter unipotent subgroups of $V$, we write the corresponding unipotent subgroup by $U_2'$ therefore
\begin{align*}
&\int\limits_{R_{\rho}^{w_0}(K)\backslash R_{\rho}(\mathbb{A})}
 \varphi(g)f_s(w_0 ug)\Tilde{\theta}_{\phi}(\tau(u)g)dudg\\
 =&\int\limits_{G_{\rho}^{w_0}(K)\backslash G_{\rho}(\mathbb{A})}\int\limits_{(U_2')^{w_0}(K)\backslash U_2'(\mathbb{A})}\int\limits_{[U_2]}
 \varphi(g)f_s(w_0 u_1u'g)\Tilde{\theta}_{\phi}(\tau(u_1)\tau(u')g)du_1du'dg\\
 =&\int\limits_{G_{\rho}^{w_0}(K)\backslash G_{\rho}(\mathbb{A})}\int\limits_{(U_2')^{w_0}(K)\backslash U_2'(\mathbb{A})}\varphi(g)f_s(w_0u'g)\of{\int\limits_{(K\backslash\mathbb{A})^2} \Tilde{\theta}_{\phi}(\tau(x_{1231}(y_3)x_{1221}(y_4))\tau(u')g)dy_3dy_4}du'dg.
\end{align*}
The last equation is because the invariant property of $f_s$ after pushing $u_1$ across $w_0$.
On the other hand, the image of $u_1\in N_P$ in $Sp_8^{\rho}$ is of the shape
$$\begin{pmatrix}
1 &   &  * & *  & * & * & * & *\\
  & 1 &  * & *  & * & * & * & *\\
  &   & 1 &   & 0 & 0 & * & *\\
  &   &   & 1 & 0 & 0 & * & *\\
  &   &   &   & 1 &   &  * & *\\
  &   &   &   &   & 1 &  * & *\\
  &   &   &   &   &   & 1 & \\
  &   &   &   &   &   &   & 1 \\
\end{pmatrix},$$
Hence by the properties of Weil representation, the most inner integral with respect to $y_3$ and $y_4$ is 
\begin{align*}
    &\int_{(K\backslash\mathbb{A})^2} \Tilde{\theta}_{\phi}(\tau(x_{1231}(y_3)x_{1221}(y_4))\tau(u')g)dy_3dy_4\\ 
    =&\sum_{\underline{\xi}\in K^4} \omega_{\psi}(\tau(u')g).\phi(\underline{\xi}) \int_{(K\backslash\mathbb{A})^2}\psi(2\xi_2 y_3-2\rho \xi_1 y_4)dy_3dy_4\\
    =&\sum_{{\xi}_3,{\xi}_4\in K} \omega_{\psi}(\tau(u')g).\phi((0,0,\xi_3,\xi_4)).
\end{align*}
Next, we write $g=n_1g_1\in G_{\rho}$ where $n_1\in N_P$ and let $P_{-}$ be the complementary subgroup of $N_P$ in $G_{\rho}$ then by the same reason as in case $w_1$ and changing the variables, here the computations are similar as in the case of $w_1$, we get
\begin{align*}
    &\int\limits_{G_{\rho}^{w_0}(K)\backslash G_{\rho}(\mathbb{A})}\int\limits_{(U_2')^{w_0}(K)\backslash U_2'(\mathbb{A})} \varphi(g)f_s(w_0u'g)\of{\int\limits_{(K\backslash\mathbb{A})^2} \Tilde{\theta}_{\phi}(\tau(x_{1231}(y_3)x_{1221}(y_4))\tau(u')g)dy_3dy_4}du'dg\\
    =&\int\limits_{P_{-}^{w_0}(K)\backslash P_{-}(\mathbb{A})}\int\limits_{(U_2')^{w_0}(K)\backslash U_2'(\mathbb{A})} \int\limits_{[N_P]}\varphi(n_1g_1)f_s(w_0u'n_1g_1)\of{\sum_{{\xi}_3,{\xi}_4\in K} \omega_{\psi}(\tau(u')n_1g_1).\phi((0,0,\xi_3,\xi_4))}dn_1du'dg_1\\
    =&\sum_{{\xi}_3,{\xi}_4\in K}\int\limits_{P_{-}^{w_0}(K)\backslash P_{-}(\mathbb{A})}\int\limits_{(U_2')^{w_0}(K)\backslash U_2'(\mathbb{A})}f_s(w_0u'g_1) \omega_{\psi}(\tau(u')g_1).\phi((0,0,\xi_3,\xi_4)\of{\int\limits_{[N_P]}\varphi(n_1g_1)dn_1}du'dg_1=0
\end{align*}
by the cuspidality of $\varphi$.

Finally, we are left with $\gamma_0=w_0w[2]$. It is not difficult to check that $R_{\rho}^{\gamma_0}$ contains the following root subgroups  
$$x_{1000}(r), \  x_{\pm0110}(r), \ x_{1110}(r), \ x_{1220}(r), \ x_{1100}(r)x_{1120}(\rho r), \ x_{1121}(r), \ x_{1231}(r).$$
Let $V_0$ be the 13 dimensional subgroup of $V$ omitting the root subgroups $x_{1121}(r)$ and $x_{1231}(r)$ and $G_1=G_{\rho}\cap \gamma_0^{-1}P\gamma_0$ and it contains the subgroup that generated by $x_{1000}(r)$, $x_{\pm0110}(r)$, $x_{1110}(r)$, $x_{1220}(r)$, $x_{1100}(r)x_{1120}(\rho r)$.

\begin{align*}
    &\int_{R_{\rho}^{\gamma_0}(K)\backslash R_{\rho}(\mathbb{A})}
 \varphi(g)f_s(\gamma_0 ug)\Tilde{\theta}_{\phi}(\tau(u)g)dudg\\
    =&\int_{G_1(K)\backslash SU_{2,2}(\mathbb{A})}\int_{V_0(\mathbb{A})}\int_{(K\backslash\mathbb{A})^2}\varphi(g)f_s(\gamma_0x_{1121}(y_2)x_{1231}(y_3)u_0g)\times\\ 
    &\hspace{8cm}\Tilde{\theta}_{\phi}(\tau(x_{1121}(y_2)x_{1231})(y_3)\tau(u_0)g)dy_2dy_3du_0dg.\\
    =&\int_{G_1(K)\backslash SU_{2,2}(\mathbb{A})}\int_{V_0(\mathbb{A})}\int_{(K\backslash\mathbb{A})^2}\varphi(g)f_s(\gamma_0u_0g) \Tilde{\theta}_{\phi}(\tau(x_{1121}(y_2)x_{1231})(y_3)\tau(u_0)g)dy_2dy_3du_0dg.
\end{align*}
By the properties of $\theta$-function, we have
\begin{align*}
    &\int_{(K\backslash\mathbb{A})^2}  \Tilde{\theta}_{\phi}(\tau(x_{1121}(y_2)x_{1231}(y_3))\tau(u_0)g)dy_2dy_3\\ 
    =&\sum_{\underline{\xi}\in K^4} \int_{(K\backslash\mathbb{A})^2}\omega_{\psi}((0,y_2,y_3,0)\tau(u_0)g).\phi(\underline{\xi})dy_2dy_3\\
    =&\sum_{\underline{\xi}\in K^4} \omega_{\psi}(\tau(u_0)g).\phi(\underline{\xi})\int_{(K\backslash\mathbb{A})^2}\psi(2\rho\xi_3 y_2+2\xi_2 y_3)dy_2dy_3\\
    =&\sum_{{\xi}_1,{\xi}_4\in K} \omega_{\psi}(\tau(u_0)g).\phi((\xi_1,0,0,\xi_4)).
\end{align*}
Back to the integral, we get

\begin{align*}
    &\int_{R_{\rho}^{\gamma_0}(K)\backslash R_{\rho}(\mathbb{A})}
 \varphi(g)f_s(\gamma_0 ug)\Tilde{\theta}_{\phi}(\tau(u)g)dudg\\
 =&\sum_{{\xi}_1,{\xi}_4\in K}\int_{G_1(K)\backslash  G_{\rho}(\mathbb{A})}\int_{V_0(\mathbb{A})}\varphi(g)f_s(\gamma_0u_0g) \omega_{\psi}(\tau(u_0)g).\phi((\xi_1,0,0,\xi_4))du_0dg.
\end{align*}
Now, the subgroup generated by $h_{1000}$, $x_{\pm 0110}$ of $G_{\rho}$, which is isomorphic to $GL_2$. Since the Levi subgroup of $Q$ acts on its unipotent radical $V$ by the adjoint action and this action induces a action of $GL_2$ on the unipotent subgroup generated by $U_{0001}$ and $U_{0111}$ with two orbits. One is 0, the other is the open orbit that contains $x_{0111}(1)$, and its isotropy subgroup in $GL_2$ is $GL_1^{\Delta}\times N_2$, where $GL_1^{\Delta}=\{h(t, t^2, t, 1)| t\in K^{\times}\}$ and $N_2=\{x_{0110}(r)|r\in K\}$. The image of $GL_1^{\Delta}$ in $SO(7)$ can be written as a matrix

$$\begin{pmatrix}
t^{-1} &   &   &   &   &   & \\
  & t^{-1} &   &   &   &   & \\
  &   & 1 &   &   &   & \\
  &   &   & 1 &   &   & \\
  &   &   &   & 1 &   & \\
  &   &   &   &   & t & \\
  &   &   &   &   &   & t\\
\end{pmatrix} \mapsto
\begin{pmatrix}
t^{-1} &   &   &   \\
  & 1 &   &   \\
  &   & 1 &  \\
  &   &   & t \\
\end{pmatrix}\subset SU_{2,2}.$$
By this we have  $$\omega_{\psi}(\tau(u_0)g).\phi((\xi_1,0,0,\xi_4))= \omega_{\psi}(g'\tau(u_0)g).\phi((0,0,0,1))$$
where $g'\in GL_2(K)$ and also by the left invariant properties $\varphi(g)=\varphi(g'g)$ and $f_s(\gamma_0u_0g))=f_s((\gamma_0g'\gamma_0^{-1})\gamma_0u_0g)$, the above integral is
\begin{align*}
   &\int_{G_1(K)\backslash G_{\rho}(\mathbb{A})}\int_{V_0(\mathbb{A})}\varphi(g)f_s(\gamma_0u_0g) \omega_{\psi}(\tau(u_0)g).\phi(0)du_0dg \ +\\
    &\sum_{g'\in (GL_1^{\Delta}\times N_2)(K)\backslash (GL_1\times SL_2)(K)}\int_{G_1(K)\backslash G_{\rho}(\mathbb{A})}\int_{V_0(\mathbb{A})}\varphi(g)f_s(\gamma_0u_0g) \omega_{\psi}(\tau(u_0)g).\phi((0,0,0,1))du_0dg
\end{align*}
For the first summand, we see $G_1$ contains an unipotent subgroup $N_1$ of $G_{\rho}$ generated by 
$$x_{1000}(r), \   \ x_{1110}(r), \ x_{1220}(r), \ x_{1100}(r)x_{1120}(\rho r).$$
Hence, we can write $g=n_1g_1\in G_{\rho}$ and $n_1\in N_1$ and then push $n_1$ to the left. By the embedding of $G_{\rho}=SU_{2,2}$ into $Sp^{\rho}_8$ above, we have $f_s(\gamma_0 u_0n_1g_1)=f_s((\gamma_0n_1\gamma_0^{-1})\gamma_0(n_1^{-1}u_0n_1)g_1)$ and $\omega_{\psi}(\tau(u_0)n_1g_1).\phi(0)=\omega_{\psi}(n_1\tau(n_1^{-1}u_0n_1)g_1).\phi(0)=\omega_{\psi}(\tau(n_1^{-1}u_0n_1)g_1).\phi(0)$ for the last two identities of Weil representation. Therefore, substituting those in to the first summand, we can change the variables on $V 
^{\gamma_0}(\mathbb{A})\backslash V(\mathbb{A})=V_0(\mathbb{A})$ and it doesn't change the Haar measure for $N_1$ is an unipotent subgroup. Now, by cuspidality, we end up with the most inner integral 
$$\int_{[N_1]}\varphi(n_1g_1)dn_1=0.$$

We are left with only the second summation. We have $G_1=(GL_2)\cdot N_1$, and let $N_0=N_1N_2$ which is the unipotent subgroup of $G_{\rho}$ as in theorem. Now, we can write $n_0\in N_0$ as
$$n_0=x_{1000}(r)x_{0110}(r_1)x_{1110}(r_2)x_{1220}(r_3)x_{1100}(r_4)x_{1120}(\rho r_4).$$ The image of $n_0$ in $Sp_8^{\rho}$ has the form
$$\begin{pmatrix}
1 &   &   &  * & * & * & * & *\\
  & 1 &  * &   & * & * & * & *\\
  &   & 1 &   & 0 & r & * & *\\
  &   &   & 1 & r & 0 & * & *\\
  &   &   &   & 1 &   &   & * \\
  &   &   &   &   & 1 &  * & \\
  &   &   &   &   &   & 1 & \\
  &   &   &   &   &   &   & 1 \\
\end{pmatrix},$$
where $r$ corresponding the unipotent subgroup $x_{1000}$. And then by the properties of Weil's representation we have
\begin{align*}
    \omega_{\psi}(\tau(u_0)n_0g).\phi((0,0,0,1))&=\omega_{\psi}(n_0)\omega_{\psi}(\tau(n_0^{-1}u_0n_0)g).\phi(0,0,0,1)\\
    &=\psi(\rho r)\omega_{\psi}(\tau(n_0^{-1}u_0n_0)g).\phi(0,0,0,1).
\end{align*}
By the fact that $\varphi$ is a cusp form $f_s$ is invariant for $K$-elements and the first formula of Weil representation. We have
\begin{align*}
    &\sum_{g'\in GL_1^{\Delta}\times N_2(K)\backslash GL_2(K)}\int_{G_1(K)\backslash G_{\rho}(\mathbb{A})}\int_{V_0(\mathbb{A})}\varphi(g'g)f_s((\gamma_0g'\gamma_0^{-1})\gamma_0u_0g)\times\\ 
    &\hspace{9cm}\omega_{\psi}(g'\tau(u_0)g).\phi((0,0,0,1))du_0dg\\
    =&\int_{GL_1^{\Delta}\times N_0(K)\backslash G_{\rho}(\mathbb{A})}\int_{V_0(\mathbb{A})}\varphi(g)f_s(\gamma_0u_0g) \omega_{\psi}(\tau(u_0)g).\phi((0,0,0,1))du_0dg\\
    =&\int_{GL_1^{\Delta}\times N_0(\mathbb{A})\backslash SU_{2,2}(\mathbb{A})}\int_{V_0(\mathbb{A})}\left(\int_{[N_0]}\varphi(n_0g)\psi(r)dn_0\right)f_s(\gamma_0u_0g)\times\\
    &\hspace{9cm}\omega_{\psi}(\tau(u_0)g).\phi((0,0,0,1))du_0dg
\end{align*}
Next, let $\psi_{\rho}(x)=\psi(\rho x)$ and $u(x)=x_{0100}(x)x_{0120}(\rho x)$ and the Fourier expansion of $\varphi(u(x)g)$ with respect to $x$, we have
$$\varphi(u(x)g)=\sum_{\beta\in K}\left(\int_{K\backslash \mathbb{A}}\varphi(u(r_1)g)\psi_{\rho}(\beta r_1)dr_1\right)\overline{\psi_{\rho}}(\beta x).$$
Letting $x=0$ and plugging it into the $N_0$ integral, we get
\begin{align*}
    &\int_{[N_0]}\varphi(n_0g)\psi_{\rho}(r)dn_0\\
    =&\sum_{\beta\in K}\int_{[N_0]}\int_{K\backslash \mathbb{A}}\varphi(u(r_1)n_0g)\psi_{\rho}(\beta r_1+r)dr_1dn_0\\
    =&\sum_{\beta\in K^*}\int_{[N_0]}\int_{K\backslash \mathbb{A}}\varphi(u(r_1)n_0g)\psi_{\rho}(\beta r_1+r)dr_1dn_0.
\end{align*}
The last equation is because if $\beta=0$, by the cuspidality we will get $0$. Since $u(r_1)$ is in the complement of $N_0$ in $N_B$ and $u(r_1)n_0$ is in the unipotent radical $N_B$, recall that $h_0(t)=h(t,t^2,t,1)$ maps to $GL_1^{\Delta}$ then making change the variable $r_1\leftrightarrow \beta r_1$ in $n_B$, we get
\begin{align*}
    &\sum_{\beta\in K^*}\int_{[N_0]}\int_{K\backslash \mathbb{A}}\varphi(u(r_1)n_0g)\psi_{\rho}(\beta r_1+r)dr_1dn_0\\
    =&\sum_{\beta\in K^*}\int_{[N_B]}\varphi(h_0(\beta^{-1})n_Bh_0(\beta)g)\psi_{\rho}(r_1+r)dn_B\\
    =&\sum_{\beta\in K^*}W_{\varphi}\left(h_0(\beta)g\right).
\end{align*}
Finally, plugging this summation into the integral, we have
\begin{align*}
    I=&\sum_{\beta\in K^*}\int_{GL_1^{\Delta}\times N_0(\mathbb{A})\backslash G_{\rho}(\mathbb{A})}\int_{V_0(\mathbb{A})}
    W_{\varphi}(h_0(\beta)g)f_s(\gamma_0u_0g)\omega_{\psi}(\tau(u_0)g).\phi((0,0,0,1))du_0dg\\
    =&\int_{N_0(\mathbb{A})\backslash G_{\rho}(\mathbb{A})}\int_{V_0(\mathbb{A})}
    W_{\varphi}(g)f_s(\gamma_0u_0g)\omega_{\psi}(\tau(u_0)g).\phi((0,0,0,1))du_0dg.
\end{align*}
The proof of the theorem is complete.

\end{proof}

It follows from the theorem \ref{main theorem-1} that $I(\varphi, \phi, f_s)$ is factorizable. More precisely, let $\pi=\otimes_{v}\pi_{v}$, $\omega_{\psi}=\otimes_{v}\omega_{\psi,v}$, $\phi=\otimes_v\phi_{v}$ and $f_s=\otimes_vf_{s,v}$ be the restricted pure tensors. If $\varphi\in \pi$ is chosen so that $W_{\varphi}=\otimes_vW_v$, then 
$$I(\varphi, \phi, f_s)=\prod_vI_{loc}(W_v, \phi_v, f_{s,v}),$$
where
$$I_{loc}(W_v, \phi_v, f_{s,v})=\int\limits_{N_0(K_v)\backslash SU_{2,2}(K_v)}\int\limits_{V_0(K_v)}W_v(g)f_{s,v}(\gamma_0 u_0g)\omega_{\psi,v}(\tau(u_0)g).\phi_{v}(0,0,0,1)du_0dg.$$

\section{The Unramified Computations}

In this section, we shall compute the unramified local integral in the non-split, quasi-split case. Now, let $K$ be a local nonarchimedean field, $\mathcal{O}_K$ its integer ring, $p$ its uniformizer and $q$ be the number of elements of residue field of $K$. We normalize all measure such that $vol(\mathcal{O_k})=1$ Since in this section we only encounter the $K$-points of algebraic groups, so we suppress "$(K)$". We assume that $\rho$ and 2 are both units in $K$ and unless stated otherwise, we always assume that $\rho$ is not a square. Let $\pi$ be an unramified irreducible generic representation of $U_{2,2}$. Let $W(g)$ be the normalized spherical vector in the Whittaker model of $\pi$, i.e. $W(gk)=W(g)$ for all $k\in U_{2,2}(\mathcal{O}_K)$ the maximal compact subgroup of $U_{2,2}$. We also notice that since $\rho$ is a unit, a $U_{2,2}$-torus conjugation identifies those Whittaker characters $\psi_{\rho}(x)$ and $\psi(x)$ so in the following we shall suppress $\rho$ and we have
$$W(x_{1000}(r_1)x_{0100}(r_2)x_{0120}(\rho r_2)n'g)=\psi(r_1+r_2)W(g),$$
where $x_{1000}(r_1)x_{0100}(r_2)x_{0120}(\rho r_2)n'$ is an element in the maximal unipotent radical of $U_{2,2}$, $\psi$ is an unramified additive character of $K$. Let $f_s$ be the spherical vector in local $\Ind_{P(K)}^{F_4(K)}\delta_P^s$. Thus, 
$$f_s(pgk)=\delta_P^s(p)f_s(g)\ \text{and} \ f_s(e)=1,$$
where $p\in P(K), g\in F_4(K)\ \text{and}\ k\in F_4(O_K)$. Finally, We also let $\phi$ be the characteristic function of $W^+(\mathcal{O})$ in $\mathcal{S}(W^+(K))$. 

We recall some fact in \cite{B-G1998} p.137. Let $B_{F_4}=TN$ be the standard Borel subgroup of $F_4$, there $T$
is the maximal torus of $F_4$ and $N$ the unipotent radical. Then we define a character $\chi_s$ of $B_{F_4}$ is given by $$\chi_s(h(t_1, t_2, t_3, t_4)n)=|t_1t_2|^{5s-1}|t_3|^{-10s+4}|t_4|^{10s-3}.$$
Let $f_s^0$ be the normalized unramified vector in $Ind_{B_{F_4}}^{F_4}\chi_s$. Moreover, for simple reflection $w\in W(F_4)$, we set 
$$(M_wf_s)(g)=\int_Kf_s(wx_{\alpha}(r)g)dr.$$
Then we have
$$M_{w[3]}M_{w[2]}: Ind_P^{F_4}\delta_P^s\to Ind_{B(F_4)}^{F_4}\chi_s.$$
Moreover,
\begin{align}\label{IntOp}
   (M_{w[3]}M_{w[2]})f_s=\frac{\zeta(5s-1)\zeta(10s-3)}{\zeta(5s)\zeta(10s-2)}f_s^0. 
\end{align}
Finally, we let $J^{\rho}(a,b)$ be the following integral
\begin{align*}
    \int\limits_{K^6}f_s^0(w_1x_{0100}(r_1)x_{1111}(r_2)x_{1222}(r_3)x_{1232}(r_4)x_{1242}(r_5)x_{1342}(r_6))\psi(a^{-1}b^2r_1+ab^{-1}r_2+\rho r_3+r_5)dr_i,
\end{align*}
where $w_1=w[412342312342]$. 

The Langlands $L$-group $^LU_{2,2}$ is defined as 
$$^LU_{2,2}=GL_4(\mathbb{C})\rtimes \Gamma_{E/K},$$
where the semidirect product action is as follow, let $^tg$ be the usual transpose of matrices. The Frobenius element $Fr$ acts on $GL_4$ as
$$Fr.g=\begin{pmatrix}
     &   &   &   1&  \\
     &   & -1  &   & \\
     &  1 &   &   & \\
   -1  &   &   &   &\\
    \end{pmatrix}{}^tg^{-1}\begin{pmatrix}
     &   &   &   -1&  \\
     &   & 1  &   & \\
     &  -1 &   &   & \\
   1  &   &   &   &\\\end{pmatrix}.$$
Now by the general theory of automorphic L-functions, see \cite{Borel}, we can associate $\pi$ with a semisimple conjugacy class and we choose a representative as $t_{\pi}=(s_{\pi}, Fr)$ in $^LU_{2,2}$. Since $s_{\pi}$ is well-defined up to Frob-conjugacy, thus we can assume that $s_{\pi}\in Sp_4(\mathbb{C})$.

Finally, we can set up the Casselman-Shalika formula for the quasisplit groups.The first discussion about such formula is due to \cite{Tamir1991}, Proposition 1, p.170, and \cite{H-G} gives a very nice interpretation of this formula. We first consider the maximal $K$-split torus $T_d$ as $Fr$-fixed points of the maximal torus $T$ of $U_{2,2}$, we shall denote this inclusion as 
$$\iota: T_d \hookrightarrow T.$$
Then it induces a surjective homomorphism on the dual side
$$\iota^*: T^{\vee}\to T_d^{\vee}$$ given by $t\mapsto\iota^*(t(Fr.t))$. In our case, since the root system of $SU_{2,2}$ with respect to the maximal $T$ is $A_3$-type. The root system $SU_{2,2}$ with respect to maximal split torus $T_d$ is of $C_2$-type and hence the root system of $^LSU_{2,2}$ with respect to the maximal torus $T_d^{\vee}$ is of $B_2$-type. So if we let
$$s_{\pi}=\begin{pmatrix}
    \mu_1 & & & \\
          & \mu_2 & & \\
          &  & \mu_2^{-1} & \\
          &  &  & \mu_1^{-1},
\end{pmatrix}\text{, then }Fr.t=\begin{pmatrix}
    \mu_1 & & & \\
          & \mu_2 & & \\
          &  & \mu_2^{-1} & \\
          &  &  & \mu_1^{-1},
\end{pmatrix},$$
therefore the image of $t$ in $T_d^{\vee}$ is 
$$\overline{s_{\pi}}=\begin{pmatrix}
    \mu_1^2 & & & & \\
          & \mu_2^2 &  & & \\
          & & 1 & &\\
         & &  & \mu_2^{-2} & \\
         & &  &  & \mu_1^{-2}\end{pmatrix}\in SO_5(\mathbb{C}).
          $$
Now, the Casselman-Shalika formula in our case takes the form
$$K(t)=\delta_B^{-\frac{1}{2}}(t_d)W(t_d)=\chi(\Gamma^{SO_5}_{m,n})(\overline{s_{\pi}}),$$
where $t_d$ corresponds to an weight element $(m,n)$ of $X_*(T_d)\cong X^*(T_d^{\vee})$ in the dominant chamber and thus gives rise to an irreducible representation $\Gamma^{SO_5}_{m,n}$ and $\chi(\Gamma^{SO_5}_{m,n})$ is the character of $\Gamma^{SO_5}_{m,n}$.

\subsection{Computation of the Zeta-Integral Side}
By Theorem \ref{main theorem-1}, at the unramified places, the local integral is 
$$I_{loc}(W,\phi, f_s)=\int_{N_0\backslash SU_{2,2}}\int_{V_0}W(g)f_s(\gamma_0 u_0g)\omega_{\psi}(\tau(u_0)g).\phi(0,0,0,1)du_0dg.$$\label{Loc int}
Then we have
\begin{proposition}\label{prop4.1}
    With the unramified data given above and $Re(s)$ large enough, we have
    \begin{align*}
        I_{loc}(W,\phi, f_s)=\frac{\zeta(5s-1)\zeta(10s-3)}{\zeta(5s)\zeta(10s-2)}\int_{T_d}K(t)J^{\rho}(a,b)\eta(a)|a|^{5s-2}dt.
    \end{align*}
Here $T_d$ is the maximal split torus of $SU_{2,2}$ and 
$t=\text{diag }(b,ab^{-1},ba^{-1}, b^{-1})\in T_d$ and $\eta=\eta_{E/K}$ is the
quadratic character of $K^{\times}$ attached to the field extension $E/K$.
\end{proposition}
\begin{proof} We will proceed in several steps.\\

\textbf{Step 1 - Iwasawa decomposition on $g$:}\\

  By the Iwasawa decomposition, we can parametrize $N_0\backslash SU_{2,2}$ as $g=x_{0100}(r)x_{0120}(\rho r)tk$, where $t$ is in the maximal non-split torus of $SU_{2,2}$ and $k$ is in the maximal compact subgroup of $SU_{2,2}$. On the other hand, if we let $T$ be the maximal torus of $SU_{2,2}$ and $T^0$ be its maximal compact subgroup, by \cite{Borel}, we have  $T/T^0\cong T_d/T_d^0$ and therefore\footnote{One can apply the same method to simplify the local computation in \cite{H2012}}, 
  \begin{align*}
      I_{loc}=&\int_{T_d}\int_{K}\int_{V_0}W(x_{0100}(r)x_{0120}(\rho r)t_d)f_s(\gamma_0 u_0x_{0100}(r)x_{0120}(\rho r)t_d)\times \\
      &\hspace{6cm}\omega_{\psi}(\tau(u_0)x_{0100}(r)x_{0120}(\rho r)t_d).\phi(0,0,0,1)du_0dr\delta_B^{-1}(t_d)dt.
  \end{align*}

  \textbf{Step 2 - Pushing $x_{0100}(r)x_{0120}(\rho r)$ to the left:}\\
  
  We have that
  $$u_0=x_{0001}(x_1)x_{0011}(x_2)x_{0121}(x_3)x_{0111}(x_4)x_{1111}(y_1)x_{1221}(y_4)u_0',$$
  here $u_0'$ is an element in $[V,V]$ i.e. the product of all 1-parameter unipotent subgroups of $x_{\alpha}$ with roots $\alpha=(***2)$. To simplify the notation, we shall denote the integrand in the $I_{loc}$ as $\Phi$. Using the Chevalley commutator relation with this notation, it will simplify the change of variables. Hence, we have 
  \begin{align*}
      &\int\limits_{V_0}\int\limits_K \Phi(u_0x_{0100}(r)x_{0120}(-\rho r))drdu_0\\
      =&\int\limits_{V_0}\int\limits_K\Phi(x_{0001}(x_1)x_{0011}(x_2)x_{0121}(x_3)x_{0111}(x_4)x_{1111}(y_1)x_{1221}(y_4)x_{0100}(r)x_{0120}(-\rho r)u_0')drdu_0\\
      =&\int\limits_{V_0}\int\limits_K\Phi(x_{0100}(r)x_{0001}(x_1)x_{0011}(x_2)x_{0121}(x_3)x_{0111}(x_4)x_{1111}(y_1)x_{0120}(\rho r)x_{1221}(y_4)u_0')drdu_0.
  \end{align*}
  We have the following commutator relation
  \begin{align*}
      x_{1111}(y_1)x_{0120}(\rho r)=x_{0120}(\rho r)x_{1231}(C_1\rho ry_1)x_{2342}(C_2\rho ry_1^2)x_{1111}(y_1),
  \end{align*}
  here the integers $C_1$ and $C_2$ are determined by the structure constants of $F_4$. Now, we can push $x_{1231}(C_1\rho ry_1)$ to the most left, $x_{2342}$ to the most right and changing the variables. So we have 
  \begin{align*}
     &\int\limits_{V_0}\int\limits_K \Phi(u_0x_{0100}(r)x_{0120}(\rho r))drdu_0\\
     =&\int\limits_{V_0}\int\limits_K\Phi (x_{0100}(r)x_{0001}(x_1)...x_{0111}(x_4)x_{0120}(\rho r)x_{1231}(C_1\rho ry_1)x_{2342}(C_2\rho ry_1^2)\\
     &\hspace{10.6cm}\times x_{1111}(y_1)x_{1221}(y_4)u_0')drdu_0\\
     =&\int\limits_{V_0}\int\limits_K\Phi(x_{1231}(C_1\rho ry_1)x_{0100}(r)x_{0001}(x_1)...x_{0111}(x_4)x_{0120}(\rho r)x_{1111}(y_1)x_{1221}(y_4)u_0'\\
     &\hspace{8.5cm} \times x_{1242}(C_1C_3\rho rx_2y_1)x_{2342}(C_2\rho ry_1^2))drdu_0\\
     =&\int\limits_{V_0}\int\limits_K\Phi(x_{1231}(C_1\rho ry_1)x_{0100}(r)x_{0120}(\rho r)x_{0001}(x_1)...x_{0111}(x_4)x_{1111}(y_1)x_{1221}(y_4)u_0')drdu_0.
  \end{align*}
  
On the other hand, we can get rid of $x_{1231}$ in $f_s$ since $\gamma_0x_{1231}\gamma_0^{-1}=x_{0011}$. Also by the properties of Weil representation, by the last integral above, we find 
  \begin{align*}
      W(x_{0100}(r)x_{0120}(\rho r)t_d)&=\psi(r)W(t_d),\\
f_s((\gamma_0x_{1231}(C_1\rho ry_1)\gamma_0^{-1})\gamma_0 x_{0100}(r)x_{0120}(\rho r) u_0t_d)&=f_s((\gamma_0x_{0120}(\rho r)\gamma_0^{-1})\gamma_0 x_{0100}(r) u_0t_d)\\
      &=f_s(x_{0122}(\rho r)\gamma_0 x_{0100}(r) u_0t_d)\\
      &=f_s(\gamma_0 x_{0100}(r) u_0t_d),\\
      \omega_{\psi}(\tau(x_{1231}(C_1\rho ry_1))x_{0100}(r)x_{0120}(\rho r)\tau(u_0)t_d).\phi(0,0,0,1)&=\omega_{\psi}(x_{0100}(r)x_{0120}(\rho r)\tau(u_0)t_d).\phi(0,0,0,1)\\
      &=\omega_{\psi}(\tau(u_0)t_d).\phi(0,0,0,1).
  \end{align*}
So our integral is reduced to
\begin{align*}
     I_{loc}=&\int_{T_d}\int_{K}\int_{V_0}W(t_d)f_s(\gamma_0 x_{0100}(r) u_0t_d)\omega_{\psi}(\tau(u_0)t_d).\phi(0,0,0,1)\psi(r)du_0dr\delta_B^{-1}(t_d)dt.
\end{align*}

\textbf{Step 3 - Pushing the torus $t_d$ to the left:}\\

The split torus element can be parametrized by cocharacters of $F_4$ as $t_d=h(a, b^2, b, 1)$, we first conjugate $t_d$ across $u_0$ and after changing the variables it produces a factor $|a|^{-1}$ from the Haar measure $du_0$, then conjugating it across $x_{0100}(r)$ we have $t_d^{-1}x_{0100}(r)t_d=x_{0100}(ab^{-2}r)$ and across $\gamma_0$ we have $\gamma_0h(a, b^2, b, 1)\gamma_0^{-1}=h(a,a,a,b)$. On the other hand, by the polarization of $W$, we have $t_d|_{W^+}=\text{diag}_E(b, ab^{-1})$ and hence ${\det}_E(t_d)=a$. Now by the properties of $f_s$ and Weil representation for unitary group, we get 
\begin{align*}
    f_s(\gamma_0 x_{0100}(r) u_0t_d)&=f_s(h(a,a,a,b)\gamma_0 x_{0100}(ab^{-2}r)u_0)\\
    &=|a|^{5s}f_s(\gamma_0 x_{0100}(ab^{-2}r)u_0),\\
    \omega_{\psi}(\tau(u_0)t_d).\phi(0,0,0,1)&=\gamma_{\psi}(a^2)\eta(a)|a|\omega_{\psi}(\tau(u_0)).\phi(0,0,0,ab^{-1})\\
    &=\gamma_{\psi}(a^2)\eta(a)|a|\phi(x_1,x_2,x_3,x_4+ab^{-1})\psi(2\rho ab^{-1}y_1+\rho r_1+r_2)\\
&=\eta(a)|a|\phi(x_1,x_2,x_3,x_4+ab^{-1})\psi(2\rho ab^{-1}y_1+\rho r_1+r_2)
\end{align*}
Back to the integral, since $2\rho$ is a unit and then changing variables and also notice that $\delta_B(t_d)=|ab^2|^2$, we have
\begin{align*}
     I_{loc}=&\int\limits_{T_d}\int\limits_{K}\int\limits_{V_0}W(t_d)|a|^{5s}f_s(\gamma_0 x_{0100}(ab^{-2}r)u_0)\eta(a)|a|\omega_{\psi}(\tau(u_0)).\phi(0,0,0,ab^{-1})\psi(r)|a|^{-1}du_0dr\delta_B^{-1}(t_d)dt\\
     =&\int\limits_{T_d}\int\limits_{K}\int\limits_{V_0}W(t_d)f_s(\gamma_0 x_{0100}(r)u_0)\phi(x_1,x_2,x_3,x_4+ab^{-1})\psi(a^{-1}b^2r+2\rho ab^{-1}y_1+\rho r_1+ r_2)\times\\
     &\hspace{11cm}\eta(a)|a|^{5s-1}|b|^2du_0dr\delta_B^{-1}(t_d)dt\\
     =&\int\limits_{T_d}\int\limits_{K}\int\limits_{V_0}K(t_d)f_s(\gamma_0 x_{0100}(r)u_0)\phi(x_1,x_2,x_3,x_4+ab^{-1})\psi(a^{-1}b^2r+2\rho ab^{-1}y_1+\rho r_1+ r_2)\times\\
     &\hspace{12.6cm}\eta(a)|a|^{5s-2}du_0drdt.\\
\end{align*}
By the Casselman-Shalika formula, if $|ab^{-1}|> 1$ we have $K_{\varphi}(t_d)=0$ hence $|x_4+ab^{-1}|\le 1$ implies $|x_4|\le 1$. Thus the integration over roots $0001$, $0011$, $0121$ and $0111$ is restricted to integers. By right $K(F_4)$-invariant property, we get
\begin{align*}
    \int_{T_d}\int_{K}\int_{V_1}K(t_d)f_s(\gamma_0 x_{0100}(r)u_0)\psi(a^{-1}b^2r+2\rho ab^{-1}y_1+\rho r_1+ r_2)\eta(a)|a|^{5s-2}du_0drdt,
\end{align*}
here $V_1$ is a 9-dimensional unipotent subgroup of $F_4$ associates roots $1111$, $1221$, $0122$, $1122$, $1222$, $1232$, $1242$, $1342$, $2342$. Now we denote the inner two integrals as
\begin{align*}
    J^{\rho}_0(a,b)=\int_{K}\int_{V_1}f_s(\gamma_0 x_{0100}(r)u_0)\psi(a^{-1}b^2r+2\rho ab^{-1}y_1+\rho r_1+r_2)du_0dr.
\end{align*}

\textbf{Step 4 - Reducing the integral over the root $x_{0122}$.}\\

Adding $x_{1100}(m)$ with $|m|\le 1$ on the right and then consider the integral
\begin{align*}
    I_1=\int_Kf_s(\gamma_0 x_{0100}(r)u_0'x_{2342}(r'))\psi(a^{-1}b^2r+2\rho ab^{-1}y_1+\rho r_1+r_2)dr'
\end{align*}
By the commutator relation
\begin{align*}
    x_{0122}(u)x_{1100}(v)=&x_{1100}(v)x_{1222}(C_{1,1}uv)x_{0122}(u)\\
    x_{1242}(u)x_{1100}(v)=&x_{1100}(v)x_{2342}(C'_{1,1}uv)x_{1242}(u)\\
    \gamma_0x_{1100}\gamma_0^{-1}=x_{1000}
\end{align*}
and then push $x_{1100}(m)$ to the left we have
$$I_1=\psi(\rho mr_{0122})I_1.$$
Thus the integral over $0122$ variable is restricted to integers. \\

\textbf{Step 5 - Reducing the integral over the root $x_{1122}$.}\\

Similarly, adding $x_{0120}(m)$ with $|m|\le 1$ on the right and pushing it to the left. The commutator relations we applied in this case are
\begin{align*}
    &\gamma_0x_{0120}\gamma_0^{-1}=x_{0122},\\
    &x_{1122}(r_1)x_{0120}(m)=x_{0120}(m)x_{1242}(mr_1)x_{1122}(r_1),\\
    &x_{1222}(r_1)x_{0120}(m)=x_{0120}(m)x_{1342}(mr_1)x_{1222}(r_1),\\
     &x_{1111}(r_1)x_{0120}(m)=x_{0120}(m)x_{1231}(mr_1)x_{2342}(mr_1^2)x_{1111}(r_1).
\end{align*}
Finally, after change notation of variables we get 
\begin{align*}
   J^{\rho}_0(a,b)=&\int\limits_{K^8}f_s(\gamma_0x_{0100}(r_1)x_{1111}(r_2)x_{1221}(r_7)x_{1222}(r_3)x_{1232}(r_4)x_{1242}(r_5)x_{1342}(r_6)x_{2342}(r_8))\times\\
    &\hspace{9cm}\psi(a^{-1}b^2r_1+ab^{-1}r_2+\rho r_3+r_5)dr_i.
\end{align*}

\textbf{Step 6 - Completing the proof.}\\

We shall write the Weyl element $\gamma_0=w[23]w_1$ and pushing $w_1$ across $x_{1221}$ and $x_{2342}$ then
\begin{align*}
   J^{\rho}_0(a,b)=&\int\limits_{K^8}f_s(w[23]w_1x_{1221}(r_7)x_{2342}(r_8)x_{0100}(r_1)x_{1111}(r_2)x_{1222}(r_3)x_{1232}(r_4)x_{1242}(r_5)x_{1342}(r_6))\times\\
    &\hspace{9cm}\psi(a^{-1}b^2r_1+2\rho ab^{-1}r_2+\rho r_3+r_5)dr_i\\
    =&\int\limits_{K^8}f_s(w[23]x_{0010}(r_7)x_{0120}(r_8)w_1x_{0100}(r_1)x_{1111}(r_2)x_{1222}(r_3)x_{1232}(r_4)x_{1242}(r_5)x_{1342}(r_6))\times\\
    &\hspace{9cm}\psi(a^{-1}b^2r_1+2\rho ab^{-1}r_2+\rho r_3+r_5)dr_i.
\end{align*}
Now by the relation of $f_s$ and $f_s^0$ (\ref{IntOp}), we get
\begin{align*}
    J^{\rho}_0&(a,b)=\frac{\zeta(5s-1)\zeta(10s-3)}{\zeta(5s)\zeta(10s-2)}J^{\rho}(a,b).
\end{align*}
Completing the proof.

\end{proof}
Next we shall compute the integral $J^{\rho}(a,b)$. We recall for $z_1, z_2, z_3\in K^*$ and $g\in F_4$, define
\begin{align*}
    F_{z_1, z_2, z_3}(g)= \int_{K^4}f_s^0(w[4121]x_{1000}(r_1)x_{0100}(r_2)x_{1100}(r_3)x_{0001}(r_4)g)\psi(z_1r_1+z_2r_2+z_3r_4)dr_i.
\end{align*}
Then the function $F_{z_1, z_2, z_3}(g)$ has the following properties, see \cite{B-G1998}, p.139.
 \begin{enumerate}[1)]
    \item $F_{z_1, z_2, z_3}(g)$ is $K(F_4)$-invariant, where $K(F_4)$ is the maximal compact subgroup of $F_4$.\label{1}
    \item For the torus $h(t_1, t_2, t_3, t_4)$ of $F_4$,
    $$F_{z_1, z_2, z_3}(h(t_1, t_2, t_3, t_4)g)=|t_1t_2|^{-5s+3}|t_3|^{20s-8}|t_4|^{-10s+5}F_{t_1^2t_2^{-1}z_1, t_1^{-1}t_2^2t_3^{-2}z_2, t_3^{-1}t_4^2z_3}(g).$$\label{2}
    \item Let $U_{F_4}$ be the maximal unipotent radical of $F_4$ and $x_{1000}(r_1)x_{0100}(r_2)x_{1100}(r_3)x_{0001}(r_4)u'\in U_{F_4}$, here $u'$ is the product of 1-unipotent subgroups of positive root that differ from $1000$, $0100$, $1100$, $0001$. Then we have
    $$F_{z_1, z_2, z_3}(x_{1000}(r_1)x_{0100}(r_2)x_{1100}(r_3)x_{0001}(r_4)u'g)=\psi^{-1}(z_1r_1+z_2r_2+z_3r_4)F_{z_1, z_2, z_3}(g).$$\label{3}
    \item $F_{z_1, z_2, z_3}(e)=0$, if $|z_i|>1$ for some $i=1,2,3$.\label{4}
    \item $F_{\epsilon_1z_1, \epsilon_2z_2, \epsilon_3z_3}(e)=F_{z_1, z_2, z_3}(e)$ for any $\epsilon_i\in \mathcal{O}_K^*$, $i=1,2,3$.\label{5}
 
 \end{enumerate}
 
 With the above $F_{z_1, z_2, z_3}(g)$ and after changing the notation of variables again, we have   
\begin{proposition}\label{prop4.2}
    Let $x=q^{-5s+2}$, we have
    \begin{align*} 
        J^{\rho}(a,b)=&\int_{K^{2}}F_{a^{-1}b^2, 1, ab^{-1}}(w[32]x_{0100}(r_1)x_{0110}(r_2))\psi(\rho r_1)dr_1dr_2\\
        =&F_{a^{-1}b^2, 1, ab^{-1}}(e)-x^3F_{p^{-1}a^{-1}b^2, 1, p^{-1}ab^{-1}}(e).
    \end{align*}
\end{proposition}
\begin{proof}
As in \cite{B-G1998}, we write $w_1=w_2w[234]$ and push $w[234]$ to the right
\begin{align*}
    J^{\rho}(a,b)=&\int_{K^6}f_s^0(w_2x_{0120}(r_1)x_{1110}(r_2)x_{1122}(r_5)x_{1242}(r_6)x_{1120}(r_3)x_{1121}(r_4))\\
    &\hspace{7cm}\times\psi(a^{-1}b^2r_1+2\rho ab^{-1}r_2+\rho r_3+r_5)dr_i,
\end{align*}
then write $w_2=w[4121]w[32413]=w[4121]w_3$ and pushing $w_3$ across first 4 roots
\begin{align*}
    J^{\rho}(a,b)=&\int_{K^6}f_s^0(w[4121]x_{1000}(r_1)x_{0001}(r_2)x_{0100}(r_5)x_{1100}(r_6)w_3x_{1120}(r_3)x_{1121}(r_4))\\
    &\hspace{7cm}\times\psi(a^{-1}b^2r_1+2\rho ab^{-1}r_2+\rho r_3+r_5)dr_i\\
    =&\int_{K^{2}}F_{a^{-1}b^2, 1, ab^{-1}}(w_3x_{1120}(r_3)x_{1121}(r_4))\psi(\rho r_3)dr_3dr_4\\
    =&\int_{K^{2}}F_{a^{-1}b^2, 1, ab^{-1}}(w[32]x_{0100}(r_3)x_{0110}(r_4))\psi(\rho r_3)dr_3dr_4
\end{align*}
This complete the proof of the first equation. Changing the indexes, we continue to compute the last integral
$$J^{\rho}(a,b)=\int_{K^{2}}F_{a^{-1}b^2, 1, ab^{-1}}(w[32]x_{0100}(r_1)x_{0110}(r_2))\psi(\rho 
r_1)dr_1dr_2.$$
We have three pieces, the first is $|r_1|\le 1$, in this case, since $\rho$ is an unit and by Iwasawa decomposition, we have
\begin{align*}
    J_1^{\rho}(a,b)=&F_{a^{-1}b^2, 1, ab^{-1}}(e)+\int_{|r_2|>1}F_{a^{-1}b^2, 1, ab^{-1}}(w[3]x_{0010}(r_2))dr_2\\
    =&F_{a^{-1}b^2, 1, ab^{-1}}(e)+\int_{|r_2|>1}F_{a^{-1}b^2, 1, ab^{-1}}(x_{0010}(r_2^{-1})h(1,1,r_2^{-1},1))dr_2\\
    =&F_{a^{-1}b^2, 1, ab^{-1}}(e)+\int_{|r_2|>1}F_{a^{-1}b^2, 1, ab^{-1}}(h(1,1,r_2^{-1},1))dr_2\\
    =&F_{a^{-1}b^2, 1, ab^{-1}}(e)+\int_{|r_2|>1}|r_2|^{-20s+8}F_{a^{-1}b^2, r_2^2, r_2ab^{-1}}(e)dr_2\\
    =&F_{a^{-1}b^2, 1, ab^{-1}}(e).
\end{align*}
The last equation is because of~\ref{4}.

The second case is $|r_1|>1$ and $|r_2|\le 1$, so by~\ref{2} and~\ref{5} we have
\begin{align*}
    J_2^{\rho}(a,b)=&\int_{|r_1|>1}F_{a^{-1}b^2, 1, ab^{-1}}(w[32]x_{0100}(r_1))\psi(\rho r_1)dr_1\\
    =&\int_{|r_1|>1}F_{a^{-1}b^2, 1, ab^{-1}}(w[3]x_{0100}(r_1^{-1})h(1,r_1^{-1},1,1))\psi(\rho r_1)dr_1\\
    =&\int_{|r_1|>1}F_{a^{-1}b^2, 1, ab^{-1}}(w[3]h(1,r_1^{-1},1,1))\psi(\rho r_1)dr_1\\
    =&\int_{|r_1|>1}F_{a^{-1}b^2, 1, ab^{-1}}(h(1,r_1^{-1},r_1^{-1},1))\psi(\rho r_1)dr_1\\
    =&\int_{|r_1|>1}|r_1|^{-15s+5}F_{r_1a^{-1}b^2, 1, r_1ab^{-1}}(e)\psi(\rho r_1)dr_1\\
    =&\sum_{k=1}^{\infty}q^{(-15s+5)k}\int_{|r_1|=q^k}F_{r_1a^{-1}b^2, 1, r_1ab^{-1}}(e)\psi(\rho r_1)dr_1\\
    =&\sum_{k=1}^{\infty}q^{(-15s+6)k}F_{p^{-k}a^{-1}b^2, 1, p^{-k}ab^{-1}}(e)\int_{|\epsilon|=1}\psi(\epsilon \rho p^{-k})d\epsilon\\
    =&-q^{-15s+5}F_{p^{-1}a^{-1}b^2, 1, p^{-1}ab^{-1}}(e).
\end{align*}
 The last equation is because if we change the variable $\rho\epsilon\leftrightarrow \epsilon$, then we get
\begin{align*}
     \int_{|\epsilon|=1}\psi(\epsilon p^{-k})d\epsilon=\begin{cases}
     -q^{-1}, & \mbox{if } k=1,\\
     0,       & \mbox{otherwise.}
     \end{cases}
\end{align*}
The last piece is $|r_1|>1$ and $|r_2|>1$.
\begin{align*}
    J_3^{\rho}(a,b)=&\int_{|r_1|>1}\int_{|r_2|>1}F_{a^{-1}b^2, 1, ab^{-1}}(w[32]x_{0100}(r_1)x_{0110}(r_2))\psi(\rho r_1)dr_1dr_2\\
    =&\sum_{k=1}^{\infty}\sum_{l=1}^{\infty}q^{k+l}\int_{|\epsilon_1|=1}\int_{|\epsilon_2|=1}F_{a^{-1}b^2, 1, ab^{-1}}(w[32]x_{0100}(\epsilon_1p^{-k})x_{0110}(\epsilon_2p^{-l}))\psi(\epsilon_1\rho p^{-k})d\epsilon_1d\epsilon_2\\
    =&\sum_{k=1}^{\infty}\sum_{l=1}^{\infty}q^{k+l}\int_{|\epsilon_1|=1}\int_{|\epsilon_2|=1}F_{a^{-1}b^2, 1, ab^{-1}}(w[3]x_{0010}(\epsilon_2p^{-l})w[2]x_{0100}(\epsilon_1p^{-k}))\times\\
    &\hspace{11cm}\psi(\epsilon_1\rho p^{-k})d\epsilon_1d\epsilon_2\\
    =&\sum_{k=1}^{\infty}\sum_{l=1}^{\infty}q^{k+l}\int_{|\epsilon_1|=1}\int_{|\epsilon_2|=1}F_{a^{-1}b^2, 1, ab^{-1}}(w[3]x_{0010}(\epsilon_2p^{-l})x_{0100}(\epsilon_1^{-1}p^{k})h(1,p^k,1,1))\times\\
    &\hspace{11cm}\psi(\epsilon_1\rho p^{-k})d\epsilon_1d\epsilon_2.
\end{align*}
We have the following identities

\begin{align*}
    &w[3]h(\epsilon_2^{-1},\epsilon_2^{-1},1,1)=h(\epsilon_2^{-1},\epsilon_2^{-1},\epsilon_2^{-1},1)\\
    &h(\epsilon_2^{-1},\epsilon_2^{-1},1,1)x_{0010}(p^{-l}) h(\epsilon_2,\epsilon_2,1,1)=x_{0010}(\epsilon_2p^{-l}) \\
    &h(\epsilon_2,\epsilon_2,1,1)x_{0100}(\epsilon_1^{-1}p^{k})h(\epsilon_2^{-1},\epsilon_2^{-1},1,1)=x_{0100}(\epsilon_2\epsilon_1^{-1}p^{k}),
    \end{align*}
so the last summation is
\begin{align*}
    &\sum_{k=1}^{\infty}\sum_{l=1}^{\infty}q^{k+l}\int_{|\epsilon_1|=1}\int_{|\epsilon_2|=1}F_{\epsilon'a^{-1}b^2, \epsilon_2, \epsilon'ab^{-1}}(w[3]x_{0010}(p^{-l})x_{0100}(\epsilon_2\epsilon_1^{-1}p^{k})h(1,p^k,1,1))\times\\
    &\hspace{12cm}\psi(\epsilon_1\rho p^{-k})d\epsilon_1d\epsilon_2.
\end{align*}
By \cite{GZ1988} we have the structure constants $N_{0010, 0100}=-1$ and $N_{0010, 0110}=-2$ and thus the commutator relation
$$x_{0010}(p^{-l})x_{0100}(\epsilon_2\epsilon_1^{-1}p^{k})=x_{0100}(\epsilon_2\epsilon_1^{-1}p^{k})x_{0110}(-\epsilon_2\epsilon_1^{-1}p^{k-l})x_{0120}(\epsilon_2\epsilon_1^{-1}p^{k-2l})x_{0010}(p^{-l}).$$
Plugging this into $F_{\epsilon'a^{-1}b^2, \epsilon_2, \epsilon'ab^{-1}}$ and by the invariant property~\ref{2} and~\ref{3} we get
\begin{align*}
    &F_{\epsilon'a^{-1}b^2, \epsilon_2, \epsilon'ab^{-1}}(w[3]x_{0010}(p^{-l})x_{0100}(\epsilon_2\epsilon_1^{-1}p^{k})h(1,p^k,1,1))\\
    =&F_{\epsilon'a^{-1}b^2, \epsilon_2, \epsilon'ab^{-1}}(x_{0100}(\epsilon_2\epsilon_1^{-1}p^{k-2l})w[3]x_{0010}(p^{-l})h(1,p^k,1,1))\\
    =&\psi^{-1}(\epsilon_2^2\epsilon_1^{-1} p^{k-2l})F_{\epsilon'a^{-1}b^2, \epsilon_2, \epsilon'ab^{-1}}(w[3]x_{0010}(p^{-l})h(1,p^k,1,1))\\
    =&\psi^{-1}(\epsilon_2^2\epsilon_1^{-1} p^{k-2l})F_{\epsilon'a^{-1}b^2, \epsilon_2, \epsilon'ab^{-1}}(h(1,p^k,p^k,1)w[3]x_{0010}(p^{k-l}))\\
    =&\psi(-\epsilon_2^2\epsilon_1^{-1} p^{k-2l})q^{(-15s+5)k}F_{p^{-k}\epsilon'a^{-1}b^2, \epsilon_2, p^{-k}\epsilon'ab^{-1}}(w[3]x_{0010}(p^{k-l})).
\end{align*}
If $k<l$, by Iwasawa decomposition we have
$$w[3]x_{0010}(p^{k-l})=x_{0010}(p^{l-k})h(1,1,p^{l-k},1)$$
and therefore by the property~\ref{2} and~\ref{4} we get $F_{p^{-k}\epsilon'a^{-1}b^2, \epsilon_2, p^{-k}\epsilon'ab^{-1}}(w[3]x_{0010}(p^{k-l}))=0$.
Now we left with $k\ge l\ge 1$ and the summation is
\begin{align*}
    &\sum_{k=1}^{\infty}\sum_{l=1}^{\infty}q^{(-15s+6)k+l}F_{p^{-k}a^{-1}b^2, 1, p^{-k}ab^{-1}}(e)\int_{|\epsilon_1|=1}\int_{|\epsilon_2|=1}\psi(-\epsilon_2^2\epsilon_1^{-1} p^{k-2l}+\epsilon_1\rho p^{-k})d\epsilon_1d\epsilon_2.
\end{align*}
For the integration over the character $\psi$ and by changing the variable $\dfrac{\epsilon_2}{\epsilon_1}\leftrightarrow  \epsilon_2$ we have
\begin{align*}
    &\int_{|\epsilon_1|=1}\int_{|\epsilon_2|=1}\psi(-\epsilon_2^2\epsilon_1^{-1}p^{k-2l}+\epsilon_1\rho p^{-k})d\epsilon_1d\epsilon_2\\
    =&\int_{|\epsilon_1|=1}\int_{|\epsilon_2|=1}\psi(\epsilon_1p^{-k}(\rho-\epsilon_2^2p^{2k-2l}))d\epsilon_2d\epsilon_1.
\end{align*}
If $k>l\ge 1$ the integral is 0 for $\rho$ is not a square and hence the integration on $\epsilon_1$ is vanishing, thus we only need to compute the case $k=l=1$. By Hensel's lemma, $\rho-\epsilon_2^2$ is an unit for $\rho$ is not a square, then we have
\begin{align*}
\int_{|\epsilon_1|=1}\int_{|\epsilon_2|=1}\psi(\epsilon_1p^{-1}(\rho-\epsilon_2^2))d\epsilon_2d\epsilon_1=\frac{q-1}{q}\int_{|\epsilon_1|=1}\psi(\epsilon_1p^{-1})d\epsilon_1=-\frac{1}{q}+\frac{1}{q^2}.
\end{align*}
And the summation equals
\begin{align*}
    -q^{-15s+6}F_{p^{-1}a^{-1}b^2, 1, p^{-1}ab^{-1}}(e)+q^{-15s+5}F_{p^{-1}a^{-1}b^2, 1, p^{-1}ab^{-1}}(e).
\end{align*}
Now, by the above computation, our unipotent integral can be written as
\begin{align*}
        \int_{K^{2}}F_{a^{-1}b^2, 1, ab^{-1}}(w[32]x_{0100}(r_1)x_{0110}(r_2))\psi(\rho r_1)dr_1dr_2
        =F_{a^{-1}b^2, 1, ab^{-1}}(e)-x^3F_{p^{-1}a^{-1}b^2, 1, p^{-1}ab^{-1}}(e).
\end{align*}

\end{proof}

\begin{remark}
    If $\rho=d^2$ is a square, the above integral over the additive character $\psi$ will be
    \begin{align*}
&\int_{|\epsilon_1|=1}\int_{|\epsilon_2|=1}\psi(\epsilon_1p^{-1}(d^2-\epsilon_2^2))d\epsilon_2d\epsilon_1\\
=&\int\limits_{\epsilon_2\equiv\pm d(\text{mod}.\mathfrak{m})}\int\limits_{|\epsilon_1|=1}d\epsilon_1d\epsilon_2+\int\limits_{\epsilon_2\not\equiv\pm d(\text{mod}.\mathfrak{m})}\int\limits_{|\epsilon_1|=1}\psi(\epsilon_1p^{-1}(d^2-\epsilon_2^2))d\epsilon_1d\epsilon_2\\
=&
\frac{2}{q}\frac{q-1}{q}+\frac{q-3}{q}(-\frac{1}{q})=\frac{q+1}{q^2}.
\end{align*}
Then the "last piece" above equals
\begin{align*}
    q^{-15s+6}F_{p^{-1}a^{-1}b^2, 1, p^{-1}ab^{-1}}(e)+q^{-15s+5}F_{p^{-1}a^{-1}b^2, 1, p^{-1}ab^{-1}}(e).
\end{align*} 
In this case, the Proposition 4.2 is reduced to one of the Bump-Ginzburg's computation in \cite{B-G1998}, p.142, identity (1).
\end{remark}

Next, we shall recall that, see \cite{B-G1998} p.144
$$H_{z_1,z_2}(g)=\int_{K^3}f_s^0(w[121]x_{1000}(r_1)x_{0100}(r_2)x_{1100}(r_3))\psi(z_1r_1+z_2r_2)dr_i.$$
Therefore,
\begin{align*}
    F_{z_1,z_2,z_3}(e)&=\int_KH_{z_1,z_2}(w[4]x_{0001}(r))\psi(z_3r)dr\\
    &=\frac{\zeta(10s-4)}{\zeta(10s-3)}(1-x^2|z_3|^{10s-4})H_{z_1,z_2}(e).
\end{align*}
We have the following
\begin{lemma}\label{lemma 4.3}
    For r and k are nonnegative integers, we have
    $$H_{p^k, p^r}(e)=\frac{\zeta(5s-2)^2\zeta(10s-4)}{\zeta(5s-1)^2\zeta(10s-3)}(1-x^{r+1})(1-x^{k+1})(1-x^{r+k+2}).$$
\end{lemma}
\begin{proof}
    See \cite{B-G1998}, p.147, Lemma 4.2.

\end{proof}
\begin{proposition}\label{prop4.4}
    Let $P_{m,n}(x)=(1-x^{n+1})[(1-x^{2m+2})(1-x^{n+2})-x^3(1-x^{2m})(1-x^{n})].$ Then the Zeta-integral has the following expansion,
    \begin{align*}
       I_{loc}&(W,\phi, f_s)=\frac{\zeta(5s-2)\zeta(10s-4)^2}{\zeta(5s)\zeta(5s-1)\zeta(10s-2)\zeta(10s-3)}\sum_{m, n\ge0}(-1)^n\chi_{SO(5)}(n,2m)x^{2m+n}P_{m,n}(x).
    \end{align*}
\end{proposition}
\begin{proof}
    First, Parametrizing the maximal split torus of $SU_{2,2}$ as $a=p^{2m+n}$ and $b=p^{m+n}$, thus $z_1=a^{-1}b^2=p^n$ and $z_3=ab^{-1}=p^m$ Let $\varpi_1$ and $\varpi_2$ be the fundamental representations of $SO_5(\mathbb{C})$ and $\varpi_1$ the 5 dimensional irreducible representation. Let $\chi_{SO(5)}(n,2m)$ be the formal character of the $n\varpi_1+2m\varpi_2$ evaluated at $\overline{s_{\pi}}$, i.e. $\chi_{SO(5)}(n,2m)=\chi(\Gamma_{n,2m}^{SO(5)})(\overline{s_{\pi}})$. The Casselman-Shalika formula gives
$$K(t)=\chi_{SO(5)}(n,2m).$$
Next, by lemma~\ref{lemma 4.3}, we have
\begin{align*}
    F_{p^n, 1, p^m}(e)=&\frac{\zeta(10s-4)}{\zeta(10s-3)}(1-x^{2m+2})H_{p^n,1}(e)\\
    =&\left(\frac{\zeta(5s-2)\zeta(10s-4)}{\zeta(5s-1)\zeta(10s-3)}\right)^2(1-x^{2m+2})(1-x)(1-x^{n+1})(1-x^{n+2})\\
    =&\frac{\zeta(5s-2)\zeta(10s-4)^2}{\zeta(5s-1)^2\zeta(10s-3)^2}(1-x^{2m+2})(1-x^{n+1})(1-x^{n+2}).
\end{align*}
Now by proposition~\ref{prop4.2} and property~\ref{4},
\begin{align*}
    J^{\rho}(a,b)=&J(p^{2m+n}, p^{m+n})=F_{p^n, 1, p^m}(e)-x^3F_{p^{n-1}, 1, p^{m-1}}(e)\\
    =&\frac{\zeta(5s-2)\zeta(10s-4)^2}{\zeta(5s-1)^2\zeta(10s-3)^2}(1-x^{n+1})[(1-x^{2m+2})(1-x^{n+2})-x^3(1-x^{2m})(1-x^{n})] 
\end{align*}
Finally, by proposition~\ref{prop4.1} and $\eta(a)=(-1)^{\ord_K(a)}=(-1)^{2m+n}=(-1)^{n}$ we have
\begin{align*}
        &I_{loc}(W,\phi, f_s)=\frac{\zeta(5s-1)\zeta(10s-3)}{\zeta(5s)\zeta(10s-2)}\int_{T_d}K_{\varphi}(t)J(a,b)\eta(a)|a|^{5s-2}dt\\
        =&\frac{\zeta(5s-2)\zeta(10s-4)^2}{\zeta(5s)\zeta(5s-1)\zeta(10s-2)\zeta(10s-3)}\\
       \times &\sum_{m, n\ge0}(-1)^n\chi_{SO(5)}(n,2m)x^{2m+n}(1-x^{n+1})[(1-x^{2m+2})(1-x^{n+2})-x^3(1-x^{2m})(1-x^{n})].
    \end{align*}
Complete the computations.

\end{proof}

\subsection{Computation of the Local Adjoint $L$-function of $SU_{2,2}$}

In this section, we shall compute the local $L$-function side. All notations are the same as before. The adjoint $L$-function of $SU_{2,2}$ is defined as follow. Let $GL_4(\mathbb{C})$ acts as conjugation on its Lie algebra $\mathfrak{sl}_4(\mathbb{C})$ and $Fr$ acts as the negative of differential of $Fr.g$, thus for $X\in \mathfrak{sl}_4(\mathbb{C})$, 
$$Fr.X=\begin{pmatrix}
     &   &   &   1&  \\
     &   & -1  &   & \\
     &  1 &   &   & \\
   -1  &   &   &   &\\\end{pmatrix}{}^tX\begin{pmatrix}
     &   &   &   -1&  \\
     &   & 1  &   & \\
     &  -1 &   &   & \\
   1  &   &   &   &\\\end{pmatrix}.$$
   It is easy to see that the action of $Fr$ is compatible with $g$. We denote this action of $GL_4(\mathbb{C})\rtimes \Gamma_{E/K}$ on $\mathfrak{sl}_4(\mathbb{C})$ as $\text{Ad}$, which we will apply to define the local adjoint $L$-function. Also, noticed that there is only one other way to define the action of $Fr$ which is compatible with the semidirect product, i.e. the differential of $Fr.g$,
   $$Fr'.X=\begin{pmatrix}
     &   &   &   1&  \\
     &   & -1  &   & \\
     &  1 &   &   & \\
   -1  &   &   &   &\\\end{pmatrix}{}^tX\begin{pmatrix}
     &   &   &   1&  \\
     &   & -1  &   & \\
     &  1 &   &   & \\
   -1  &   &   &   &\\\end{pmatrix},$$
and in this case, as was discussed in \cite{H2012}, we should inserting the quadratic character associated to the field extension $E/K$ into the induction data for the Eisenstein series and we denote this action by $\text{Ad}\otimes\chi_{\Gamma_{E/K}}$. The local adjoint $L$-function is given by
$$L(s,\pi,Ad)=\dfrac{1}{det(I_{15}-Ad(t_{\pi})q^{-s})}.$$
Also, let
$$\zeta(s)=\frac{1}{1-q^{-s}}$$
be the local Zeta-function.

Now, if we consider $Sp_4(\mathbb{C})\times \Gamma_{E/K}$ acts on $\mathfrak{sl}_4(\mathbb{C})$ by the restriction of the Ad of $^LU_{2,2}$ above, we have 
$$\mathfrak{sl}_4(\mathbb{C})=(\Gamma_{2,0}^{\Sp_4}\boxtimes\chi_{\Gamma_{E/K}})\oplus(\Gamma_{0,1}^{\Sp_4}\boxtimes\mathbb{I}_{\Gamma_{E/K}}),$$
where $\Gamma_{2,0}^{\Sp_4}=\text{Sym}^2(\Gamma_{1,0}^{\Sp_4})$ is the 10-dimensional adjoint representation of $Sp_4(\mathbb{C})$. Then we have, on the first piece, $Fr$ acts with eigenvalue $-1$ and $s_{\pi}$ acts with eigenvalues 
$$1,1, \mu_1^2, \mu_1^{-2}, \mu_1\mu_2, (\mu_1\mu_2)^{-1}, \mu_1\mu_2^{-1}, \mu_1^{-1}\mu_2, \mu_2^2, \mu_2^{-2},$$ 
and on the second piece, $Fr$ acts with eigenvalue $1$ and $s_{\pi}$ acts with eigenvalues 
$$1, \mu_1\mu_2, (\mu_1\mu_2)^{-1}, \mu_1\mu_2^{-1}, \mu_1^{-1}\mu_2.$$ 
Now by the definition of local $L$-functions and recall that $x=q^{-5s+2}$, the local adjoint $L$-function is
\begin{align*}
   &L(5s-2, \pi, Ad)=\frac{1}{(1+x)(1-x^2)}\\
   &\times\frac{1}{(1+\mu_1^2x)(1+\mu_1^{-2}x)(1-\mu_1^2\mu_2^2x^2)(1-\mu_1^{-2}\mu_2^{-2}x^2)(1-\mu_1^{2}\mu_2^{-2}x^2)(1-\mu_1^{-2}\mu_2^2x^2)(1+\mu_2^2x)(1+\mu_2^{-2}x)}. 
\end{align*}
We rewrite the above $L$-function as
\begin{align}\label{(1)}
    (1+x)L(5s-2,\pi,Ad)=&\frac{1}{(1+\mu_1^2x)(1+\mu_1^{-2}x)(1+\mu_2^2x)(1+\mu_2^{-2}x)}\times \nonumber\\ 
    &\frac{1}{(1-\mu_1^2\mu_2^2x^2)(1-\mu_1^{-2}\mu_2^{-2}x^2)(1-x^2)(1-\mu_1^{-2}\mu_2^2x^2)(1-\mu_1^2\mu_2^{-2}x^2)}.
\end{align}
Let $\Pi$ be the unramified representation of $SO_5$ whose semisimple conjugacy class contains 
$$t_{\Pi}=\text{diag}(\mu_1^2,\mu_2^2,\mu_2^{-2},\mu_1^{-2})\in Sp_4(\mathbb{C}).$$
Then, as before, let $\Gamma_{1,0}^{Sp_4}$ and $\Gamma_{0,1}^{Sp_4}$ be the standard 4-dimensional representation and the 5-dimensional irreducible representation of $Sp_4(\mathbb{C})$ respectively. We recall that the representation $\Gamma_{0,1}^{Sp_4}$ can be constructed as follow, $\bigwedge^2\Gamma_{1,0}^{Sp_4}\cong \mathbb{C}\oplus\Gamma_{0,1}^{Sp_4}$, see \cite{Fulton1991RepresentationTA}, \S{16.2}. We shall use $\varpi_1$ and $\varpi_2$ to denote the corresponding weights, respectively. Thus for a general representation $\Gamma_{m,n}^{\Sp_4}$, its weight $(m,n)=m\varpi_1+n\varpi_2$. For a finite dimensional complex representation $r$ we write $\chi(r)$ for the character of $r$ and thus $\chi(r)(s_r)$ is the trace of an element $r(s_r)$.

Now, by (\ref{(1)}) we can write
$$L(5s-2, \Pi, \Gamma_{1,0}^{Sp_4})=\frac{1}{(1+\mu_1^2x)(1+\mu_1^{-2}x)(1+\mu_2^2x)(1+\mu_2^{-2}x)},$$
and
$$L(10s-4, \Pi, \Gamma_{0,1}^{Sp_4})=\frac{1}{(1-\mu_1^2\mu_2^2x^2)(1-\mu_1^{-2}\mu_2^{-2}x^2)(1-x^2)(1-\mu_1^{-2}\mu_2^2x^2)(1-\mu_1^2\mu_2^{-2}x^2)}.$$
On the other hand, by the definition of local Langlands $L$-functions, we have the right hand side of the above are
\begin{align*}
    L(5s-2, \Pi, \Gamma_{1,0}^{Sp_4})=\sum_{i=0}^{\infty}\chi(\text{Sym}^i(\Gamma_{1,0}^{Sp_4}))(t_{\Pi})(-1)^ix^i,
\end{align*}
and
\begin{align*}
    L(10s-4, \Pi, \Gamma_{0,1}^{Sp_4})=\sum_{j=0}^{\infty}\chi(\text{Sym}^j(\Gamma_{0,1}^{Sp_4}))(t_{\Pi})x^{2j}.
\end{align*}
Therefore, we need to study the decomposition of the symmetric powers in those coefficients. In fact, we have
\begin{proposition}
    $$\text{Sym}^k(\Gamma_{1,0}^{Sp_4})\cong\Gamma_{k,0}^{Sp_4}$$ for all integers $k\ge 0$.
\end{proposition}
\begin{proof}
    This is easy. See \cite{Fulton1991RepresentationTA}, \S 16.2, Exercise 16.10, p.249.

\end{proof}

\begin{proposition}
    $$\text{Sym}^k(\Gamma_{0,1}^{Sp_4})\cong\bigoplus_{\substack{j=0\\ j\equiv k \pmod{2}}}^k\Gamma_{0,j}^{Sp_4}$$ for all integers $k\ge 0$.
\end{proposition}
\begin{proof}
    This is because for any integer $k\ge 2$ we have a surjective map $\text{Sym}^k(\Gamma_{0,1}^{Sp_4})\to \text{Sym}^{k-2}(\Gamma_{0,1}^{Sp_4})$ and its kernel is $\Gamma_{0,k}^{Sp_4}$, then applying the induction on $k$, see \cite{Fulton1991RepresentationTA}, \S 16.2, Exercise 16.11, p.250.

\end{proof}
By the above two propositions, we have 
\begin{align}\label{(2)}
    L(5s-2, \Pi, \Gamma_{1,0}^{Sp_4})=&\sum_{i=0}^{\infty}\chi(\text{Sym}^i(\Gamma_{1,0}^{Sp_4}))(t_{\Pi})(-1)^ix^i\nonumber\\
    =&\sum_{i=0}^{\infty}\chi(\Gamma_{i,0}^{Sp_4})(t_{\Pi})(-1)^ix^i,
\end{align}
and
\begin{align}\label{(3)}
    L(10s-4, \Pi, \Gamma_{0,1}^{Sp_4})=&\sum_{j=0}^{\infty}\chi(\text{Sym}^j(\Gamma_{0,1}^{Sp_4}))(t_{\Pi})x^{2j}\nonumber\\
    =&\sum_{t_1=0}^{\infty}\big(\sum_{k_1=0}^{t_1}\chi(\Gamma_{0,2k_1}^{Sp_4})(t_{\Pi})\big)x^{4t_1}+\sum_{t_2=0}^{\infty}\big(\sum_{k_2=0}^{t_2}\chi(\Gamma_{0,2k_2+1}^{Sp_4})(t_{\Pi})\big)x^{4t_2+2}\nonumber \\
    =&\sum_{k_1=0}^{\infty}\chi(\Gamma_{0,2k_1}^{Sp_4})(t_{\Pi})\big(\sum_{t_1=k_1}^{\infty}x^{4t_1}\big)+\sum_{k_2=0}^{\infty}\chi(\Gamma_{0,2k_2+1}^{Sp_4})(t_{\Pi})\big(\sum_{t_2=k_2}^{\infty}x^{4t_2+2} \big)\nonumber\\
    =&\frac{1}{1-x^4}\sum_{k_1=0}^{\infty}\chi(\Gamma_{0,2k_1}^{Sp_4})(t_{\Pi})x^{4k_1}+\frac{1}{1-x^4}\sum_{k_2=0}^{\infty}\chi(\Gamma_{0,2k_2+1}^{Sp_4})(t_{\Pi})x^{4k_2+2}\nonumber\\
    =&\frac{1}{1-x^4}\sum_{j=0}^{\infty}\chi(\Gamma_{0,j}^{Sp_4})(t_{\Pi})x^{2j}.
\end{align}
Putting (\ref{(2)}) and (\ref{(3)}) into (\ref{(1)}) and changing the notation, we get
\begin{align}\label{(4)}
    (1+x)(1-x^4)L(5s-2,\pi, Ad)=&\sum_{k=0}^{\infty}\sum_{l=0}^{\infty}\chi(\Gamma_{k,0}^{Sp_4})(t_{\Pi})\chi(\Gamma_{0,l}^{Sp_4})(t_{\Pi})(-1)^kx^{k+2l}\nonumber\\
    =&\sum_{k=0}^{\infty}\sum_{l=0}^{\infty}\chi(\Gamma_{k,0}^{Sp_4}\otimes\Gamma_{0,l}^{Sp_4})(t_{\Pi})(-1)^{k}x^{k+2l}.
\end{align}
Therefore, in order to get a closed-form expression in terms of characters of $Sp_4(\mathbb{C})$ we also need to decompose $\Gamma_{i,0}^{Sp_4}\otimes\Gamma_{0,j}^{Sp_4}$ into irreducibles. In fact, we have the following
\begin{align}\label{10}
    \chi(\Gamma_{k,0}^{Sp_4}\otimes\Gamma_{0,l}^{Sp_4})=\sum_{s=0}^{\min\{k,l\}}\sum_{i=0}^{\min\{k-s,s\}}\chi(\Gamma_{k-2i,l-s+i}^{Sp_4}).
\end{align}
We will proof this in the last section, see corollary (\ref{1.4}).

Next express the local adjoint $L$-function in terms of the characters of irreducible representations of $SO(5)$. Let $\varepsilon_1$ and $\varepsilon_2$ be the standard basis element in the weight lattice $\Lambda$ of $SO(5)$ or $Sp_4(\mathbb{C})$. In other word, we can write any weight as $[\lambda_1,\lambda_2]=\lambda_1\varepsilon_1+\lambda_2\varepsilon_2$, where $\lambda_1$ and $\lambda_2$ are integers in case $Sp_4$ and all integers or all half-integers in case $SO(5)$. Hence we have $\rho_{Sp_4}=2\varepsilon_1+\varepsilon_2$ and $\rho_{SO(5)}=\dfrac{3}{2}\varepsilon_1+\dfrac{1}{2}\varepsilon_2$. We have the following relation between the character of $SO(5)$ and $Sp_4(\mathbb{C})$. 
\begin{proposition}\label{2.1}
    With the above notations, we have
    \begin{align}\label{(6)}
        \chi_{Sp_4}([\lambda_1,\lambda_2])=\frac{\chi_{SO(5)}([\lambda_1+\dfrac{1}{2},\lambda_2+\dfrac{1}{2}])}{\chi_{SO(5)}([\dfrac{1}{2},\dfrac{1}{2}])}.
    \end{align}
\end{proposition}
\begin{proof}
    Let $s(T)=T-T^{-1}$ and by Weyl character formula, see \cite{Fulton1991RepresentationTA}, \S 24, p.400.
    \begin{align*}
                \chi_{Sp_4}([\lambda_1,\lambda_2])=\dfrac{A_{\lambda+\rho_{Sp_4}}}{A_{\rho_{Sp_4}}}.
    \end{align*}
    On the other hand, for arbitrary weight $\mu=[\mu_1,\mu_2]$ of $Sp_4(\mathbb{C})$ we have
    $$A_{\mu}=\begin{vmatrix}
        s(x_1^{\mu_1+2}) & s(x_2^{\mu_1+2})\\
        s(x_1^{\mu_2+1}) & s(x_2^{\mu_2+1})
    \end{vmatrix}=\begin{vmatrix}
        s(x_1^{(\mu_1+\frac{1}{2})+\frac{3}{2}}) & s(x_2^{(\mu_1+\frac{1}{2})+\frac{3}{2}})\\
        s(x_1^{(\mu_2+\frac{1}{2})+\frac{1}{2}}) & s(x_2^{(\mu_2+\frac{1}{2})+\frac{1}{2}})
    \end{vmatrix}.$$
    Applying the above identity and after dividing both numerator and denominator of $\chi_{Sp_4}([\lambda_1,\lambda_2])$ by $A_{\rho_{SO(5)}}$, we get the desired formula.
    
\end{proof}

Returning to the $L$-function, we shall write $\chi_{\Sp_4}(i,j)=\chi(\Gamma_{i,j}^{Sp_4})(t_{\Pi})$ to simplify the notation. By Proposition~\ref{2.1} this is the same as evaluating the right hand side of (\ref{(6)}) at 
$$\overline{s_{\pi}}=\begin{pmatrix}
    \mu_1^2 & & & & \\
     & \mu_2^2 & & & \\
     & & 1 & & \\
     & & & \mu_2^{-2} & \\
     & & & & \mu_1^{-2}
\end{pmatrix}.$$
Then we have
\begin{proposition}
     \label{4.14}The local adjoint $L$-function has the following power series expansion in terms of $SO(5)$ characters 
    \begin{align}\label{Lseries}
    &L(5s-2,\pi, Ad)=\frac{1}{(1+x)(1-x^2)^2(1-x^4)}\sum_{r,t=0}^{\infty}\dfrac{\chi_{SO(5)}(r,2t+1)(\overline{s_{\pi}})}{\chi_{SO(5)}(0,1)(\overline{s_{\pi}})}(-1)^rx^{r+2t}(1-x^{2r+2})(1-x^{2t+2}).
    \end{align}
\end{proposition}

\begin{proof}
    First, by (\ref{(4)}), (\ref{10}) and Proposition \ref{2.1} we have
    \begin{align*}
    (1+x)(1-x^4)&L(5s-2,\pi, Ad)=\sum_{k=0}^{\infty}\sum_{l=0}^{\infty}\chi(\Gamma_{k,0}^{Sp_4}\otimes\Gamma_{0,l}^{Sp_4})(t_{\Pi})(-1)^kx^{k+2l}\\
    =&\sum_{k=0}^{\infty}\sum_{l=0}^{\infty}\Big(\sum_{s=0}^{\min\{k,l\}}\sum_{i=0}^{\min\{k-s,s\}}\chi_{\Sp_4}(k-2i,l-s+i)\Big)(-1)^kx^{k+2l}\\
    =&\sum_{k=0}^{\infty}\sum_{l=0}^{\infty}\Big(\sum_{s=0}^{\min\{k,l\}}\sum_{i=0}^{\min\{k-s,s\}}\chi_{\Sp_4}([k+l-s-i,l-s+i])\Big)(-1)^kx^{k+2l}\\
    =&\sum_{k=0}^{\infty}\sum_{l=0}^{\infty}\Big(\sum_{s=0}^{\min\{k,l\}}\sum_{i=0}^{\min\{k-s,s\}}\dfrac{\chi_{SO(5)}([k+l-s-i+\dfrac{1}{2},l-s+i+\dfrac{1}{2}])(\overline{s_{\pi}})}{\chi_{SO(5)}([\dfrac{1}{2},\dfrac{1}{2}])(\overline{s_{\pi}})}\Big)(-1)^kx^{k+2l}\\
    =&\sum_{k=0}^{\infty}\sum_{l=0}^{\infty}\sum_{s=0}^{\min\{k,l\}}\sum_{i=0}^{\min\{k-s,s\}}\dfrac{\chi_{SO(5)}(k-2i,2(l-s+i)+1)(\overline{s_{\pi}})}{\chi_{SO(5)}(0,1)(\overline{s_{\pi}})}(-1)^kx^{k+2l}.
\end{align*}
Now, we denote $$A(k-2i,l-s+i)=\dfrac{\chi_{SO(5)}(k-2i,2(l-s+i)+1)(\overline{s_{\pi}})}{\chi_{SO(5)}(0,1)(\overline{s_{\pi}})}.$$
For a fixed $k$, consider the middle two summations, then
$$\sum_{l=0}^{\infty}\sum_{s=0}^{\min\{k,l\}}=\sum_{s=0}^{k}\sum_{l=s}^{\infty}.$$
Therefore the above summation reduce to
\begin{align*}
    &\sum_{k=0}^{\infty}\sum_{l=0}^{\infty}\sum_{s=0}^{\min\{k,l\}}\sum_{i=0}^{\min\{k-s,s\}}A(k-2i,l-s+i)(-1)^kx^{k+2l}\\
    =&\sum_{k=0}^{\infty}\sum_{s=0}^{k}\sum_{l=s}^{\infty}\sum_{i=0}^{\min\{k-s,s\}}A(k-2i,l-s+i)(-1)^kx^{k+2l}\\
    =&\sum_{s=0}^{\infty}\sum_{k=s}^{\infty}\sum_{l=s}^{\infty}\sum_{i=0}^{\min\{k-s,s\}}A(k-2i,l-s+i)(-1)^kx^{k+2l}.\\
\end{align*}
Then, let $k=k_1+s$ and $l=l_1+s$, and changing the range of the summations, we have
\begin{align*}
    &\sum_{s=0}^{\infty}\sum_{k=s}^{\infty}\sum_{l=s}^{\infty}\sum_{i=0}^{\min\{k-s,s\}}A(k-2i,l-s+i)(-1)^kx^{k+2l}\\
    =&\sum_{s=0}^{\infty}\sum_{k_1=0}^{\infty}\sum_{l_1=0}^{\infty}\sum_{i=0}^{\min\{k_1,s\}}A(k_1+s-2i,l_1+i)(-1)^{k_1+s}x^{k_1+2l_1+3s}\\
    =&\sum_{s=0}^{\infty}\sum_{l_1=0}^{\infty}\sum_{i=0}^{s}\sum_{k_1=i}^{\infty}A(k_1+s-2i,l_1+i)(-1)^{k_1+s}x^{k_1+2l_1+3s}\\
    =&\sum_{i=0}^{\infty}\sum_{l_1=0}^{\infty}\sum_{s=i}^{\infty}\sum_{k_1=i}^{\infty}A(k_1+s-2i,l_1+i)(-1)^{k_1+s}x^{k_1+2l_1+3s}
\end{align*}
Let $k_1=k_2+i$ and $s=s_1+i$, then the last summation is
\begin{align*}
\sum_{i=0}^{\infty}\sum_{l_1=0}^{\infty}\sum_{s_1=0}^{\infty}\sum_{k_2=0}^{\infty}A(k_2+s_1,l_1+i)(-1)^{k_2+s_1}x^{k_2+2l_1+3s_1+4i}.
\end{align*}
Finally, let $r=k_2+s_1$ and $t=l_1+i$, 
\begin{align*}
&\sum_{i=0}^{\infty}\sum_{l_1=0}^{\infty}\sum_{s_1=0}^{\infty}\sum_{k_2=0}^{\infty}A(k_2+s_1,l_1+i)(-1)^{k_2+s_1}x^{k_2+2l_1+3s_1+4i}\\
=&\sum_{r=0}^{\infty}\sum_{t=0}^{\infty}A(r,t)(-1)^rx^{r+2t}\sum_{s_1=0}^{r}\sum_{i=0}^{t}x^{2s_1+2i}\\
=&\dfrac{1}{(1-x^2)^2}\sum_{r=0}^{\infty}\sum_{t=0}^{\infty}A(r,t)(-1)^rx^{r+2t}(1-x^{2r+2})(1-x^{2t+2}).
\end{align*}
Hence, the proposition.

\end{proof}

\subsection{Comparison}
In this section we will compare the Zeta-integral and the local adjoint $L$-function. Recall that
$$P_{m,n}(x)=(1-x^{n+1})[(1-x^{2m+2})(1-x^{n+2})-x^3(1-x^{2m})(1-x^{n})].$$ 
Then by the Proposition~\ref{prop4.4} and after multiplying the normalizing factor of Eisenstein series as in \cite{B-G1998} we have
\begin{align}
        \zeta(5s)\zeta(5s-1)&\zeta(10s-2)\zeta(10s-3)\zeta(10s-4)\zeta(15s-6)\zeta(20s-8)I_{loc}(W,\phi, f_s)\nonumber\\
        &=\zeta(5s-2)\zeta(10s-4)^3\zeta(15s-6)\zeta(20s-8)\sum_{m, n\ge0}(-1)^n\chi_{SO(5)}(n,2m)x^{2m+n}P_{m,n}(x)\nonumber\\
        &=\frac{1}{(1-x)(1-x^4)(1-x^3)(1-x^2)^3}\sum_{m, n\ge0}(-1)^n\chi_{SO(5)}(n,2m)x^{2m+n}P_{m,n}(x).\label{S1}
\end{align}

By the proposition~\ref{4.14}, to reduce the $L$-function side, we need to compute $\dfrac{\chi_{SO(5)}(l,2k+1)}{\chi_{SO(5)}(0,1)}$ for any weight $(l,2k+1)$ of $SO(5)$. In fact, we have the following
\begin{align}\label{14-1}
        \dfrac{\chi_{SO(5)}(l,2k+1)}{\chi_{SO(5)}(0,1)}=\sum_{i=0}^{k}\sum_{j=0}^l(-1)^{l+k-i-j}\chi_{SO(5)}(k-i+j,2i).
    \end{align}
We will proof it in the last section, see Proposition~\ref{4.13}. According to (\ref{14-1}), we have
\begin{align*}
    &\sum_{r=0}^{\infty}\sum_{t=0}^{\infty}\dfrac{\chi_{SO(5)}(r,2t+1)(\overline{s_{\pi}})}{\chi_{SO(5)}(0,1)(\overline{s_{\pi}})}(-1)^rx^{r+2t}(1-x^{2r+2})(1-x^{2t+2})\\
    =&\sum_{r=0}^{\infty}\sum_{t=0}^{\infty}\sum_{j=0}^{r}\sum_{i=0}^{t}(-1)^{2r+t-i-j}\chi_{SO(5)}(t-i+j,2i)(\overline{s_{\pi}})x^{r+2t}(1-x^{2r+2})(1-x^{2t+2})\\
    =&\sum_{j=0}^{\infty}\sum_{i=0}^{\infty}\sum_{r=j}^{\infty}\sum_{t=i}^{\infty}(-1)^{t-i-j}\chi_{SO(5)}(t-i+j,2i)(\overline{s_{\pi}})x^{r+2t}(1-x^{2r+2})(1-x^{2t+2}).
\end{align*}
Let $r=j+r_1$ and $t=i+t_1$ and then $j_1=j+t_1$, so the last one above equals
\begin{align*}
    &\sum_{j=0}^{\infty}\sum_{i=0}^{\infty}\sum_{r_1=0}^{\infty}\sum_{t_1=0}^{\infty}(-1)^{t_1+j}\chi_{SO(5)}(j+t_1,2i)(\overline{s_{\pi}})x^{j+2i+r_1+2t_1}(1-x^{2j+2+2r_1})(1-x^{2i+2+2t_1})\\
=&\sum_{i=0}^{\infty}\sum_{r_1=0}^{\infty}\sum_{t_1=0}^{\infty}\sum_{j_1=t_1}^{\infty}(-1)^{j_1}\chi_{SO(5)}(j_1,2i)(\overline{s_{\pi}})x^{j_1+2i+r_1+t_1}(1-x^{2j_1+2+2r_1-2t_1})(1-x^{2i+2+2t_1})\\
    =&\sum_{i=0}^{\infty}\sum_{j_1=0}^{\infty}(-1)^{j_1}\chi_{SO(5)}(j_1,2i)(\overline{s_{\pi}})x^{j_1+2i}\Big(\sum_{r_1=0}^{\infty}\sum_{t_1=0}^{j_1}x^{r_1+t_1}(1-x^{2j_1+2+2r_1-2t_1})(1-x^{2i+2+2t_1})\Big).
\end{align*}
To simplify the notation, denote $m=i, n=j_1, r=r_1, t=t_1$, 
$$\sum_{m=0}^{\infty}\sum_{n=0}^{\infty}(-1)^n\chi_{SO(5)}(n,2m)(\overline{s_{\pi}})x^{2m+n}\Big(\sum_{r=0}^{\infty}\sum_{t=0}^{n}x^{r+t}(1-x^{2n+2+2r-2t})(1-x^{2m+2+2t})\Big)$$
The power series in the parentheses can be reduced as follow,
\begin{align*}
    &\sum_{r=0}^{\infty}\sum_{t=0}^{n}x^{r+t}(1-x^{2n+2+2r-2t})(1-x^{2m+2+2t})\\
    =&\sum_{t=0}^{n}x^{t}(1-x^{2m+2+2t})\sum_{r=0}^{\infty}x^{r}(1-x^{2n+2+2r-2t})\\
    =&\sum_{t=0}^{n}x^{t}(1-x^{2m+2+2t})\Big(\frac{1}{1-x}-\frac{x^{2n+2-2t}}{1-x^3}\Big)\\
    =&\frac{1}{1-x^3}\sum_{t=0}^{n}x^{t}(1-x^{2m+2+2t})(1+x+x^2-x^{2n+2-2t}).
\end{align*}
Now, we compute the above summation 
\begin{align*}
    &\sum_{t=0}^{n}x^{t}(1-x^{2m+2+2t})(1+x+x^2-x^{2n+2-2t})\\
    =&\sum_{t=0}^{n}(x^t+x^{t+1}+x^{t+2}-x^{2n+2-t})-\sum_{t=0}^{n}(x^{2m+2+3t}+x^{2m+3+3t}+x^{2m+4+3t})+\sum_{t=0}^{n} x^{2m+2n+4+t}\\
    =&\frac{(1-x^{n+1})(1+x+x^2-x^{n+2})}{1-x}-\frac{x^{2m+2}(1-x^{3n+3})}{1-x}+\frac{x^{2m+2n+4}(1-x^{n+1})}{1-x}\\
    =&\frac{(1-x^{n+1})(1+x+x^2-x^{n+2}-x^{2m+2}-x^{2m+n+3})}{1-x}=\frac{P_{m,n}(x)}{(1-x)^2}
\end{align*}
Therefore put all these together, we get 
\begin{align*}
     &\sum_{r=0}^{\infty}\sum_{t=0}^{\infty}\dfrac{\chi_{SO(5)}(r,2t+1)(\overline{s_{\pi}})}{\chi_{SO(5)}(0,1)(\overline{s_{\pi}})}(-1)^rx^{r+2t}(1-x^{2r+2})(1-x^{2t+2})\nonumber\\
     =&\frac{1}{(1-x)^2(1-x^3)}\sum_{m=0}^{\infty}\sum_{n=0}^{\infty}(-1)^n\chi_{SO(5)}(n,2m)(\overline{s_{\pi}})x^{2m+n}P_{m,n}(x)
\end{align*}
Now by Proposition \ref{4.14} we have 
\begin{align}\label{16}
     L(5s-2,\pi, Ad)=\frac{1}{(1-x)(1-x^2)^3(1-x^3)(1-x^4)}\sum_{m,n=0}^{\infty}(-1)^n\chi_{SO(5)}(n,2m)(\overline{s_{\pi}})x^{2m+n}P_{m,n}(x)
\end{align}
Then by (\ref{S1}) and (\ref{16}) we have 
\begin{theorem}\label{mainTheorem-2}
For all unramified data, we have
   \begin{align*}
      I_{loc}(W,\phi, f_s)=\frac{L(5s-2,\pi, Ad)}{\zeta(5s)\zeta(5s-1)\zeta(10s-2)\zeta(10s-3)\zeta(10s-4)\zeta(15s-6)\zeta(20s-8)}.
   \end{align*}
\end{theorem}
\renewcommand{\thetheorem}{4.10B}
\begin{theorem}
     If one consider $f'_s\in \Ind_{P(K)}^{F_4(K)}((\eta\circ\varpi_2)\delta_P^s)$ and let $x=\eta(p)q^{-5s+2}=-q^{-5s+2}.$
The local $L$-function denotes by $L(\eta,s)=(1-\eta(p)q^{-s})^{-1}$.
Then by the same computation, we have similar identity
    \begin{align*}
      I_{loc}&(W,\phi, f'_s)=\\
      &\frac{L(5s-2,\pi, \text{Ad}\otimes\chi_{\Gamma_{E/K}})}{L(\eta,5s)L(\eta,5s-1)L(\eta^2,10s-2)L(\eta^2,10s-3)L(\eta^2,10s-4)L(\eta^3,15s-6)L(\eta^4,20s-8)}.
   \end{align*}
\end{theorem}
\begin{proof}
    This amounts to change the variable $x$ to $-x$ both in the identity of Proposition~\ref{prop4.4} and the definition of local adjoint $L$-function $L(s, \pi, Ad)$ and then applying the comparison computation.

\end{proof}

\subsection{Decomposition of Symmetric Algebras}

In this section, we shall prove the identities (\ref{10}) and (\ref{14-1}).

We start with (\ref{10}). To do this, we would like to compute the Newell-Littlewood coefficients for $Sp_6(\mathbb{C})$. We now recall some facts from \cite{KOIKE1987466} that we will apply to get this decomposition. The Newell-Littlewood numbers are defined as
\begin{align*}
N_{\lambda,\mu}^{\nu}=\sum_{\alpha,\beta,\gamma}\mathrm{LR}_{\alpha,\beta}^{\lambda}\mathrm{LR}_{\alpha,\gamma}^{\mu}\mathrm{LR}_{\beta,\gamma}^{\nu}
\end{align*}
where the indices are in the set of all partitions $\mathcal{P}=\{[\lambda_1,\lambda_2, ... ,\lambda_n]\in \mathbb{Z}^n: \lambda_1\ge\lambda_2\ge ... \ge\lambda_n\ge 0\}.$
Here, $\mathrm{LR}_{\alpha,\beta}^{\lambda}$ is the Littlewood-Richardson coefficient, see \cite{macdonald1998symmetric}, chapter I, p.143. If we let $\chi_{Sp}(\lambda)$ be the universal character of symplectic group associated a partition $\lambda$, as it was shown in \cite{KOIKE1987466} \S 2.3, p.501, we have 
$$\chi_{Sp}(\lambda)\chi_{Sp}(\mu)=\sum_{\nu}N_{\lambda,\mu}^{\nu}\chi_{Sp}(\nu).$$
Now, if one apply the so-called "specialization homomorphism" defined in \cite{KOIKE1987466}, \S 2.2, p.495. As a corollary, we have if two partitions $\lambda$ and $\mu$ are in the stable range of $Sp_{2n}(\mathbb{C})$, i.e. the sum of the depth $d(\lambda)+d(\mu)\le n$, then
\begin{align*}
    \chi_{Sp_{2n}}(\lambda)\chi_{Sp_{2n}}(\mu)=\sum_{\nu}N_{\lambda,\mu}^{\nu}\chi_{Sp_{2n}}(\nu).
\end{align*}
For our purpose, we shall apply this identity for the partitions $\lambda=[k]$ and $\mu=[l,l]$. In terms of weights of $Sp_6(\mathbb{C})$ we have $\lambda$ is $(k,0,0)$ and $\mu$ is $(0,l,0)$. We have
\begin{proposition}
   Let $\Gamma_{k,0,0}$ and $\Gamma_{0,l,0}$ be the irreducible representations of $Sp_6(\mathbb{C})$ corresponding to the partitions $\lambda$ and $\mu$ respectively. Then
   \begin{align}\label{(5)}
       \chi(\Gamma_{k,0,0})\chi(\Gamma_{0,l,0})=\sum_{s=0}^{\min\{k,l\}}\sum_{i=0}^{\min\{l,k-s\}}\sum_{j=\max\{0,i-s\}}^{\min\{i,l-s\}}\chi(\Gamma_{k-2i+j,l-s+i-2j,j}).
   \end{align}
\end{proposition}
\begin{proof}
    We notice that for partitions $\lambda$ and $\mu$ we have $d(\lambda)+d(\mu)=3$ so we can use the above identity. Next, we will compute $N_{\lambda,\mu}^{\nu}$. First, assume $k\ge l$, By the definition, $N_{\lambda,\mu}^{\nu}$ is non-zero if and only if each Littlewood-Richardson coefficient is non-zero. Consider $\mathrm{LR}_{\alpha,\beta}^{\lambda}$, where the Young diagram of partition $\lambda$ consists of one row with $k$ boxes:
$$\lambda: \underbrace{\begin{ytableau}
{} & {} & \dots & {} &{}
\end{ytableau}}_{k\text{ boxes}}$$
If $\mathrm{LR}_{\alpha,\beta}^{\lambda}\neq 0$ then we must have $\alpha\subset \lambda$ and $\beta\subset \lambda$ this means the Young diagram of $\alpha$ and $\beta$ are 
$$\alpha: \underbrace{\begin{ytableau}
{} & {} & \dots & {} &{}
\end{ytableau}}_{s\text{ boxes}}$$
where $0\le s\le k$ and the skew tableaux $T$ of the shape $\lambda-\alpha$ should have the weight $\beta$ such that its word $w(T)$ is a lattice permutation (\cite{macdonald1998symmetric}, p.143), thus there is only one way to write such $T$
$$T: \underbrace{\begin{ytableau}
1 & 1 & \dots & 1 &1
\end{ytableau}}_{k-s\text{ boxes}}$$
Now we get all partitions $\beta=[k-s]$ such that $\mathrm{LR}_{\alpha,\beta}^{\lambda}=1$.

Then consider $\mathrm{LR}_{\alpha,\gamma}^{\mu}$. The Young diagram of $\mu=[l,l]$ consists of two rows and each row has $l$ boxes
$$\mu: \underbrace{\begin{ytableau}
{} & {} & \dots & {} &{}\\
{} & {} & \dots & {} &{}\\
\end{ytableau}}_{l\text{ boxes}}$$
Since we assumed that $k\ge l$ then all the possible $\alpha\subset \mu$ are satisfying $0\le s \le l$ and thus for each such $\alpha$ we have a skew-tableaux $T_{\alpha}$ with shape $\mu-\alpha$, now if $\mathrm{LR}_{\alpha,\gamma}^{\mu}\neq 0$ we must have the weight of $T_{\alpha}$ is $\gamma$ and its word $w(T_{\alpha})$ is again a lattice permutation. Therefore we have for each given $\alpha$ the only possible $T_{\alpha}$ is 
$$T_{\alpha}: \underbrace{\begin{ytableau}
\none & \none & \none & 1 & \dots & 1\\
1 & \dots & 1 & 2 & \dots & 2\\
\end{ytableau}}_{l\text{ boxes}}$$
where in the first row has $l-s$ 1's. Hence we get all partitions $\gamma=[l,l-s]$ such that $\mathrm{LR}_{\alpha,\gamma}^{\mu}=1$.

Finally, we compute $\mathrm{LR}_{\beta,\gamma}^{\nu}$ where $\beta=[k-s]$ and $\gamma=[l,l-s]$. Now the possible Young diagram of $\nu$ has 3 rows. For each such $\nu$, by the definition of Littlewood-Richardson coefficients, we are going to count the number skew tableaux $T$ of the shape $\nu-\beta$ that has weight $\gamma$ and the word $w(T)$ is a lattice permutation.
$$\nu: \begin{ytableau}
*(gray!50) & *(gray!50) & *(gray!50) & *(gray!50)\dots & *(gray!50) &*(gray!50) &1 & 1 & 1 &\dots & 1\\
1 & 1 &\dots & 1 & 2 & \dots &2 &\dots & 2\\
2 & \dots & 2
\end{ytableau}$$
where the gray boxes represent the Young diagram of $\beta$ and let $i$ and $j$ denote respectively the number of boxes labeled 1 added to the second row and the number of boxes labeled 2 added to the third row. These two parameters determine all the possible of $v$ since if $T$ is a Young tableaux and the word $w(T)$ is a lattice permutation we must have 
$$0\le i\le \min\{l,k-s\}$$ and
$$\max\{0, i-s\}\le j\le \min\{i, l-s\}.$$
The left part of the second inequality is because the number of 2s in the second row cannot greater than the number of 1s in the first row. So if $l-s>l-i$ there will be at least $i-s$ 2s going to the third row and the shape of $\nu$ is determined uniquely by $T$. Hence for each possible $\nu$, we have $\mathrm{LR}_{\beta,\gamma}^{\nu}=1$ and we have proved if $k\ge l$, 
\begin{align*}
       \chi(\Gamma_{k,0,0})\chi(\Gamma_{0,l,0})=\sum_{s=0}^l\sum_{i=0}^{\min\{l,k-s\}}\sum_{j=\max\{0,i-s\}}^{\min\{i,l-s\}}\chi(\Gamma_{k-2i+j,l-s+i-2j,j}).
   \end{align*}
If $k\le l$ by the similar argument as above one can show
\begin{align*}
       \chi(\Gamma_{k,0,0})\chi(\Gamma_{0,l,0})=\sum_{s=0}^k\sum_{i=0}^{k-s}\sum_{j=\max\{0,i-s\}}^{i}\chi(\Gamma_{k-2i+j,l-s+i-2j,j}).
   \end{align*}
Completing the proof.

\end{proof}
Now, if we specialize to group $Sp_4(\mathbb{C})$, i.e. from $\chi_{Sp_6}(x_1,x_2,x_3)$ to $\chi_{Sp_4}(x_1,x_2)$, then only the term with $j=0$ in (\ref{(5)}) contributes as a character of $Sp_4(\mathbb{C})$. We get

\begin{corollary}\label{1.4}
     Let $\Gamma_{k,0}^{Sp_4}$ and $\Gamma_{0,l}^{Sp_4}$ be the irreducible representations of $Sp_4(\mathbb{C})$ corresponding to the partitions $\lambda$ and $\mu$ respectively. Then
   \begin{align*}
       \chi(\Gamma_{k,0}^{Sp_4})\chi(\Gamma_{0,l}^{Sp_4})=\sum_{s=0}^{\min\{k,l\}}\sum_{i=0}^{\min\{k-s,s\}}\chi(\Gamma_{k-2i,l-s+i}^{Sp_4}).
   \end{align*}
\end{corollary}

Next, we proof the identity (\ref{14-1}). First, we have the following decomposition.
\begin{proposition}
    Let $\chi_{\Sp_6}(i,j,k)$ be the formal character of $Sp_6(\mathbb{C})$ corresponding the weight $(i,j,k)$. If $k,l$ are non-negative integers, then
    \begin{align*}
        \chi_{\Sp_6}(2k,l,0)\chi_{\Sp_6}(1,0,0)=&\chi_{\Sp_6}(2k+1,l,0)+\chi_{\Sp_6}(2k-1,l+1,0)\nonumber\\
        &+\chi_{\Sp_6}(2k,l-1,1)+\chi_{\Sp_6}(2k-1,l,0)+\chi_{\Sp_6}(2k+1,l-1,0).
    \end{align*}
\end{proposition}
\begin{proof}
    The same idea as Proposition 4.12, we would like to compute $N_{\lambda,\mu}^{\nu}$ for $\lambda$ and $\mu$ being $[2k+l,l]$ and $[1]$ respectively. Now by the similar argument as Proposition 4.12, we can find all possible partitions $\beta$ and $\gamma$ such that $\text{LR}_{\alpha,\beta}^{\lambda}$ and $\text{LR}_{\alpha,\gamma}^{\mu}$ are all non-zero. Here the only nontrivial case is when the Young diagram of $\alpha$ consists one box, then there are two ways to fill in the skew-Young diagram $\lambda-\alpha$ by numbers such that the corresponding skew-tableaux satisfying the word is a lattice permutation, i.e.
$$T_1: \begin{ytableau}
\none & 1 & 1 & 1 &\dots & 1\\
2 & 2 &\dots & 2\\
\end{ytableau}$$
where the first row contains $2k+l-1$ 1s and the second row contains $l$ 2s, and 
$$T_2: \begin{ytableau}
\none & 1 & 1 & 1 &\dots & 1\\
1 & 2 &\dots & 2\\
\end{ytableau}$$
where the second row contains $l-1$ 2s. Finally, the non-zero $\text{LR}_{\beta,\gamma}^{\nu}$ are given by the following partitions: For $\beta=\lambda$ and $\gamma=[1]$, we get $\nu$ are $[2k+l+1,l]$, $[2k+l,l+1]$, $[2k+l,l,1]$. For $\beta=[2k+l-1,l]$ or $\beta=[2k+l,l-1]$ and $\gamma=[\emptyset]$, we get $\nu=[2k+l-1,l]$ or $[2k+l,l-1]$. The proof is complete.

\end{proof}
If we specialize to $Sp_4(\mathbb{C})$ and applying the dual between type $B_2$ and type $C_2$, i.e. there is a $1-1$ corresponding between weights $(m,n)_{SO(5)}$ and $(n,m)_{Sp_4}$. Then we get
\begin{corollary}\label{4.16}
    For the group $Sp_4(\mathbb{C})$ and $SO(5)$ we have
    \begin{align*}
        \chi_{\Sp_4}(2k,l)&\chi_{\Sp_4}(1,0)\\
        &=\chi_{\Sp_4}(2k+1,l)+\chi_{\Sp_4}(2k-1,l+1)+\chi_{\Sp_4}(2k-1,l)+\chi_{\Sp_4}(2k+1,l-1); 
    \end{align*}
    and
    \begin{align*}
        \chi_{SO(5)}&(l,2k)\chi_{SO(5)}(0,1)\nonumber\\
        &=\chi_{SO(5)}(l,2k+1)+\chi_{SO(5)}(l+1,2k-1)+\chi_{SO(5)}(l,2k-1)+\chi_{\SO(5)}(l-1,2k+1),
    \end{align*}
where $k,l$ are all non-negative integers.
\end{corollary}

Finally, let $k$ run through all integers from $0$ to $k$ and multiply each identity by $(-1)^{k-i}$, then we have  
\begin{align*}
    \sum_{i=0}^{k}(-1)^{k-i}\chi_{SO(5)}(l+k-i,2i)=\dfrac{\chi_{SO(5)}(l,2k+1)}{\chi_{SO(5)}(0,1)}+\dfrac{\chi_{SO(5)}(l-1,2k+1)}{\chi_{SO(5)}(0,1)}.
\end{align*}
If $l=0$ then
\begin{align*}
    \dfrac{\chi_{SO(5)}(0,2k+1)}{\chi_{SO(5)}(0,1)}=\sum_{i=0}^{k}(-1)^{k-i}\chi_{SO(5)}(k-i,2i).
\end{align*}
And, therefore, we get
\begin{proposition}\label{4.13}
    Let $k\ge 0$ be an integer, then
    \begin{align*}
        \dfrac{\chi_{SO(5)}(l,2k+1)}{\chi_{SO(5)}(0,1)}=\sum_{i=0}^{k}\sum_{j=0}^l(-1)^{l+k-i-j}\chi_{SO(5)}(k-i+j,2i).
    \end{align*}
\end{proposition}

\section*{Acknowledgments}
The author would like to express his sincere gratitude to Dr. Hundley for his invaluable guidance, encouragement, and continuous support throughout the course of this work. The author is also grateful to Dr. Liu for carefully checking the detailed computations. The author would further like to thank Dr. Friedberg for his helpful comments and suggestions, and Dr. Ginzburg for many insightful discussions.

\section*{AI Disclosure}
Claude was used to assist with literature searches, English language editing, and checking mathematical computations.

\bibliographystyle{amsalpha}
\bibliography{mybib}

\vspace{5ex}

\noindent Department of Mathematics, Purdue University, West Lafayette, IN 47907, USA.

\noindent E-mail address: {\tt ma1104@purdue.edu}

\end{document}